\documentclass[reqno]{amsart}
\usepackage{hyperref}

\begin{document}

\title[\hfilneg \hfil  Some results in the uniqueness of meromorphic function    ]
{ Some results in the uniqueness of meromorphic function  }

\author[Xiaohuang Huang \hfil \hfilneg]
{Xiaohuang Huang}

\thanks{Email:huangxiaohuang2023@email.szu.edu.cn }

\begin{abstract}
In this research notes, we investigate some remain problems in the uniqueness of meromorphic function. Using some deep results of Yamanoii, we obtain some results in this notes.
\end{abstract}

\maketitle
\numberwithin{equation}{section}
\newtheorem{theorem}{Theorem}[section]
\newtheorem{lemma}[theorem]{Lemma}
\newtheorem{remark}[theorem]{Remark}
\newtheorem{corollary}[theorem]{Corollary}
\newtheorem{example}[theorem]{Example}
\newtheorem{problem}[theorem]{Problem}
\allowdisplaybreaks

\maketitle

\section{Some Notations }
Throughout this paper, by meromorphic functions we will always mean meromorphic functions in the complex plane. We adopt the standard notations in the
Nevanlinna theory of meromorphic functions as explained in  \cite{h3,l,y2}.
\par
\vskip 2mm Let  $f$ be a meromorphic function on the whole complex plane. For $0<r<R$, we define the following functions.\\
$$m(r,f)=\frac{1}{2\pi}\int_{0}^{2\pi}\log^{+}\left\vert f(re^{i\theta}) \right\vert\, {\rm d}\theta$$
is the average of the positive logarithmic $\left\vert f \right\vert$ on the circle $\left\vert z \right\vert=r$.\\
$$N(r,f)=\int_{0}^{r}\frac{n(t,f)-n(0,f)}{t}\, {\rm d}t+n(0,f)\log r,$$
is called the counting function of poles of $f$, where $n(t,f)$ denotes the number of poles of $f(z)$ on the disc $\left\vert z \right\vert\leq t$, multiple  poles are
counted according to their multiplicities. $n(0,f)$ denotes the multiplicity of pole of $f$ at the origin (if $f(0) \neq \infty$, then $n(0, f) = 0$). We denote $\overline{N}(r,f)$  the reduced counting function of $f$ whose the multiplicity of poles  only counts once.\\
$$T(r,f)=m(r,f)+N(r,f).$$
\par
\vskip 2mm We call $T(r, f)$  the characteristic function of $f$, and  it is easy to see that $T(r, f)$ is a non-negative function.\\
\par
\vskip 2mm Let $a$ be a complex number. Obviously, $\frac{1}{f(z)-a}$
is meromorphic on the disc $|z| \leq R$. Similar to above definitions, we   define the following functions.\\
$$m\left(r,\frac{1}{f-a}\right)=\frac{1}{2\pi}\int_{0}^{2\pi}\frac{1}{\log^{+}\left\vert f(re^{i\theta}) -\alpha\right\vert}\, {\rm d}\theta$$
is the average of the positive logarithmic $\left\vert \frac{1}{f-a} \right\vert$ on the circle $\left\vert z \right\vert=r$.\\
$$N\left(r,\frac{1}{f-a}\right)=\int_{0}^{r}\frac{n\left(t,\frac{1}{f-a}\right)-n\left(0,\frac{1}{f-\alpha}\right)}{t}\, {\rm d}t+n(0,f)\log r,$$
where $n\left(t,\frac{1}{f-a}\right)$ denotes the number of zeros of $f(z)-a$ on the disc $\left\vert z \right\vert\leq t$, multiple  poles are counted according to their multiplicities. $n\left(0,\frac{1}{f-a}\right)$ denotes the multiplicity of zeros of $f-\alpha$ at the origin. We denote $\overline{N}\left(r,\frac{1}{f-a}\right)$  the reduced counting function of $f-a$ whose the multiplicity of zero  only counts once.\\
$$T\left(r,\frac{1}{f-a}\right)=m\left(r,\frac{1}{f-a}\right)+N\left(r,\frac{1}{f-a}\right).$$
$T\left(r,\frac{1}{f-a}\right)$ is  the characteristic function of $\frac{1}{f-a}$.
\par
\vskip 2mm  For a non-constant meromorphic function $f$, we denote  by $S(r,f)$ any quantity satisfying $S(r, f) = o(T(r,f))$,  as   $r\not\in E\subset(0,\infty)$ and $r\to\infty$, where $E$ is a set of finite linear measure.  Let $f$ and $\beta$ be two meromorphic functions in the complex plane. We say that $b$ is a small function of $f$, if $T(r,b)=o(T(r,f))$, as $r\not\in E$ and $r\to\infty$. We say that two non-constant meromorphic functions $f$ and $g$ share small function $b$ CM(IM), if $f-b$ and $g-b$ have the same zeros counting multiplicities (ignoring multiplicities).  In addition, we say that $f$ and $g$ share $\infty$ CM, if $1/f$ and $1/g$ share $0$ CM, and we say that $f$ and $g$ share $\infty$ IM, if $1/f$ and $1/g$ share $0$ IM,  this can be found in \cite{y1}. We say that two non-constant meromorphic functions $f$ and $g$  share small function $b$  IM* if
$$\overline{N}\left(r,\frac{1}{f-b}\right)+\overline{N}\left(r,\frac{1}{g-b}\right)-2\overline{N}_{0}\left(r,b;f,g\right)=S(r,f)+S(r,g),$$
where $\overline{N}_{0}\left(r,b;f,g\right)$ denotes the reduced counting function of the common zeros of $f-b$ and $g-b$ in $|z|<r$. We say that two non-constant meromorphic functions $f$ and $g$  share small function $b$  CM* if
$$N\left(r,\frac{1}{f-b}\right)+N\left(r,\frac{1}{g-b}\right)-2N_{0}\left(r,b;f,g\right)=S(r,f)+S(r,g),$$
where $N_{0}\left(r,b;f,g\right)$ denotes the  counting function of the common zeros of $f-b$ and $g-b$ in $|z|<r$.
\par
\vskip 2mm We define the spherical characteristic function $\hat{T}(r,f)$ by
$$\hat{T}(r,f)=\frac{1}{\pi}\int_{1}^{r}A(t,f)\frac{\rm dt}{t},$$
where
$$A(t,f)=\int_{\mathbb{C}(t)}f^{*}\omega_{\hat{\mathbb{C}}}.$$
Here
$$\omega_{\hat{\mathbb{C}}}=\frac{1}{(1+|\omega|^{2})^{2}}\frac{\sqrt{-1}}{2}d\omega\wedge d\bar{\omega}$$
is the spherical area form on the Riemann sphere $\hat{\mathbb{C}}$ such that the total area of the Riemann sphere is $\pi$.
\par
\vskip 2mm Let $a\in\widehat{\mathbb{C}}$, we denote the counting function $\hat{N}(r,a,f)$ by
$$\hat{N}(r,a,f)=\int_{1}^{r}\hat{n}(t,a,f)\frac{\rm dt}{t},$$
where $\hat{n}(t, a, f)$ is the number of solutions to $f(z) = a$ on $\mathbb{C}(t)$ counting multiplicity. We also define the reduced counting function $\overline{N}(r,a,f)$ by
$$\overline{\hat{N}}(r,a,f)=\int_{1}^{r}\overline{\hat{n}}(t,a,f)\frac{\rm dt}{t},$$
where $\overline{\hat{n}}(t, a, f)$ is the number of solutions to $f(z) = a$ on $\mathbb{C}(t)$  without counting multiplicity.
\par
\vskip 2mm We define the chordal distance between two points in the complex plane by
$$[a,b]=\frac{|a-b|}{\sqrt{1+a^{2}}\sqrt{1+b^{2}}}.$$
We extend the chordal distance continuously by
$$[a,\infty]=\frac{1}{\sqrt{1+a^{2}}}.$$
We define the proximity function $\hat{m}(r,a,f)$ by
$$\hat{m}(r,a,f)=\int_{0}^{2\pi}\log\frac{1}{[f(re^{i\theta}),a]}\frac{\rm d\theta}{2\pi}.$$
By the first main theorem, we have
$$\hat{T}(r,f)=\hat{m}(r,a,f)+\hat{N}(r,a,f)-\hat{m}(1,a,f).$$
\par
\vskip 2mm  From Shimizu-Ahlfors characteristic [\cite{go}, Proof of Theorem 4.2, Page 20] we have
\begin{align}
&\hat{m}(r,a,f)=m\left(r,\frac{1}{f-a}\right)+S(r,f)\notag\\
&\hat{N}(r,a,f)=N\left(r,\frac{1}{f-a}\right)+S(r,f)\notag\\
&\hat{T}(r,f)=T(r,f)+S(r,f).\label{eq1}
\end{align}
\par
\vskip 2mm  For any distinct complex values $a_{1}, a_{2}, a_{3},\ldots,a_{q}$, K. Yamanoi\cite{ya} introduced the following modification of the proximity function
$$\bar{m}_{0,q}(r,f)=\sup_{(a_{1},a_{2},\ldots,a_{q})\in(\widehat{\mathbb{C}})^{q}}\int_{0}^{2\pi}\max_{1\leq j\leq q}\log\frac{1}{[f(re^{i\theta}),a_{j}]}\frac{\rm d\theta}{2\pi}.$$
\par
\vskip 2mm On the one hand, we can see from  [\cite{y4}, Remark in Page 76] that
\begin{align}
 \sum_{1\leq j\leq q}\hat{m}(r,a_{j},f)=\int_{0}^{2\pi}\max_{1\leq j\leq q}\log\frac{1}{[f(re^{i\theta}),a_{j}]}\frac{d\theta}{2\pi}+O(1)\leq \bar{m}_{0,q}(r,f).\label{eq2}
\end{align}
\par
\vskip 2mm On the other hand, it is easy to see from above, [\cite{ya}, (1.10)-(1.11), Proof of Theorem 1.6, Page 708] and \eqref{eq1} that
\begin{align}
&\bar{m}_{0,q}(r,f)=\int_{0}^{2\pi}\max_{1\leq j\leq q}\log\frac{1}{[f(re^{i\theta}),a_{j}]}\frac{\rm d\theta}{2\pi}+O(1)\notag\\
&\leq \sum_{1\leq j\leq q}\hat{m}(r,a_{j},f)=\sum_{1\leq j\leq q}m\left(r,\frac{1}{f-a_{j}}\right)+S(r,f).\label{eq3}
\end{align} 
\par
\vskip 2mm Moreover, let $a$ be a finite complex number, we denote $N_{1}\left(r,\frac{1}{f-a}\right)$, $N_{1}\left(r,\frac{1}{f-\infty}\right)$ and $\hat{N}_{1}(r,a,f)$ by

\begin{align}
&N_{1}\left(r,\frac{1}{f-a}\right)=N\left(r,\frac{1}{f-a}\right)-\overline{N}\left(r,\frac{1}{f-a}\right)\notag\\
&N_{1}\left(r,\frac{1}{f-\infty}\right)=N_{1}(r,f)=N(r,f)-\overline{N}(r,f)\notag\\
&\hat{N}_{1}(r,a,f)=N_{1}\left(r,\frac{1}{f-a}\right)+S(r,f).\label{eq4}
\end{align}

\section{Meromorphic functions that share four small functions}
\par
\vskip 2mmWe will touch only a few aspects of the uniqueness of meromorphic functions which share values or small functions with another meromorphic function we need in the following. 
\par
\vskip 2mmIn 1926, Nevanlinna \cite{n1} proved the following two celebrated  theorems.
\par
\vskip 2mm{\bf Theorem A.}\, Let $f$ and $g$ be two non-constant meromorphic functions, and let $a_{1}, a_{2}, a_{3}, a_{4}, a_{5}$ be five distinct constants. If $f$ and $g$ share $a_{1}, a_{2}, a_{3}, a_{4}, a_{5}$ IM, then $f\equiv g$.
\par
\vskip 2mm{\bf Theorem B.}\, Let $f$ and $g$ be two non-constant meromorphic functions, and let $a_{1}, a_{2}, a_{3}, a_{4}$ be four distinct constants. If $f$ and $g$ share $a_{1}, a_{2}, a_{3}, a_{4}$ CM. Then $f$ is  a M$\ddot{o}$bius transformation of $g$, two of the shared values, say $a_{1}$, $a_{2}$ are Picard values, and the cross ratio $(a_{1},a_{2},a_{3},a_{4})=-1$.
\par
\vskip 2mm Example of  Theorem B is   that we can consider two entire functions $e^{z}$ and $e^{-z}$. It is easy to check $e^{z}$ and $e^{-z}$ share $0,\infty,1,-1$ CM.
\par
\vskip 2mm Li-Qiao \cite{lq} improved five values to five small functions. They proved the following result.
\par
\vskip 2mm{\bf Theorem A1.}\, Let $f$ and $g$ be two non-constant meromorphic functions, and let $a_{1}, a_{2}, a_{3}, a_{4}, a_{5}$ be five distinct small functions of $f$ and $g$. If $f$ and $g$ share $a_{1}, a_{2}, a_{3}, a_{4}, a_{5}$ IM, then $f\equiv g$.
\par
\vskip 2mm  K. Ishizaki \cite{i} proved a weak version of  Theorem A1. 
\par
\vskip 2mm{\bf Theorem A2.}\, Let $f$ and $g$ be two non-constant meromorphic functions, and let $a_{1}, a_{2}, a_{3}, a_{4}, a_{5}$ be five distinct small functions of $f$ and $g$. If $f$ and $g$ share $a_{1}, a_{2}, a_{3}, a_{4}, a_{5}$ IM, and if $\overline{N}\left(r,\frac{1}{f-a_{5}}\right)\neq S(r,f)$, then $f\equiv g$.
\par
\vskip 2mm In 1979, G. Gundersen \cite{g1} improved  Theorem B from "4CM" to "3CM+1IM".
\par
\vskip 2mm{\bf Theorem C.}\, Let $f$ and $g$ be two non-constant meromorphic functions, and let $a_{1}, a_{2}, a_{3}, a_{4}$ be four distinct constants. If $f$ and $g$ share $a_{1},a_{2}, a_{3}$  CM, and $a_{4} $ IM, then $f$ and $g$ share $a_{i}(i=1,2,3,4)$ CM.
\par
\vskip 2mm G. Gundersen \cite{g1} gave the following example, which shows that there exists two distinct meromorphic functions that sharing four distinct values IM, but there is no relationship between them.
\par
\vskip 2mm {\bf Example} \cite{g1}\quad Let $h$ be a non-constant entire function, $a\neq0$ a finite constant. We set $f=\frac{e^{h}+a}{(e^{h}-a)^{2}}$ and $g=\frac{(e^{h}+a)^{2}}{8a^{2}(e^{h}-a)}$. By computation, we can verify that $f$ and $g$ share four values $0,\infty, \frac{1}{a},-\frac{1}{8a}$ IM. But $f$ is not a M$\ddot{o}$bius transformation of $g$.
\par
\vskip 2mm Li-Yang \cite{ly} proved the following  results, which generalize  Theorem B and  Theorem C to small functions.
\par
\vskip 2mm{\bf Theorem B1.}\, Let $f$ and $g$ be two non-constant meromorphic functions, and let $a_{1}, a_{2}, a_{3}, a_{4}$ be four distinct small functions of $f$ and $g$. If $f$ and $g$ share $a_{1}, a_{2}, a_{3}, a_{4}$ CM, then $f$ is  a M$\ddot{o}$bius transformation of $g$.
\par
\vskip 2mm{\bf Theorem C1.}\, Let $f$ and $g$ be two non-constant meromorphic functions, and let $a_{1}, a_{2}, a_{3}, a_{4}$ be four distinct small functions of $f$ and $g$. If $f$ and $g$ share $a_{1},a_{2}, a_{3}$  CM* and $a_{4} $ IM*, then $f$ and $g$ share $a_{1}, a_{2}, a_{3}, a_{4}$ CM.
\par
\vskip 2mm In 1983,  Gundersen \cite{g2} generalized  Theorem B from "3CM+1IM" to "2CM+2IM".
\par
\vskip 2mm{\bf Theorem D.}\, Let $f$ and $g$ be two non-constant meromorphic functions, and let $a_{1}, a_{2}, a_{3}, a_{4}$ be four distinct constants. If $f$ and $g$ share $a_{1}$, $a_{2}$  IM and $a_{3}, a_{4}$ CM, then $f$ and $g$ share $a_{1}, a_{2}, a_{3}, a_{4}$ CM.
\par
\vskip 2mm Steinmetz \cite{n2} indicated that there exist examples of three distinct non-constant meromorphic functions that share four values where none of the four values are shared CM for any two of the three functions.
\par
\vskip 2mm And then, Li \cite{l} and Yao \cite{y4} improved  Theorem D from four distinct values to four distinct small functions. 
\par
\vskip 2mm{\bf Theorem D1.}\, Let $f$ and $g$ be two non-constant meromorphic functions, and let $a\neq0,1,\infty$ be a  small function of $f$ and $g$. If $f$ and $g$ share $0$, $1$  IM* and $a, \infty$ CM*, then $f$ and $g$ share $0,1,a,\infty$ CM, and $f$ and $g$ satisfy one of the following seven relations: \\
\quad (i)  $f\equiv g$.\\
\quad (ii)  $f\equiv -g$ holds if $a\equiv-1$, and $1, -1$ are Picard exceptional values of $f$ and $g$.\\
\quad (iii)  $f\equiv 2-g$ holds if $a\equiv2$, and $0,2$ are Picard exceptional values of $f$ and $g$.\\
\quad (iv)  $f\equiv 1-g$ holds if $a\equiv\frac{1}{2}$, and $0,1$ are Picard exceptional values of $f$ and $g$.\\
\quad (v)  $\left(f-\frac{1}{2}\right)\left(g-\frac{1}{2}\right)\equiv 1$ holds if $a\equiv\frac{1}{2}$, and $\frac{1}{2},\infty$ are Picard exceptional values of $f$ and $g$.\\
\quad (vi)  $(f-1)(g-1)\equiv 1$ holds if $a\equiv2$, and $1,\infty$ are Picard exceptional values of $f$ and $g$.\\
\quad (vii)  $fg\equiv 1$ holds if $a\equiv-1$, and $0,\infty$ are Picard exceptional values of $f$ and $g$.\\
\par
\vskip 2mm So, according to the above theorems and examples, it is natural to raise the following conjecture.
\par
\vskip 2mm{\bf Conjecture 1.1} If two non-constant meromorphic functions share three values IM and share a fourth value CM, then do the functions necessarily share all four values CM?
\par
\vskip 2mm In 1989, E. Muse  \cite{m1} introduced the following definition.
\par
\vskip 2mm{\bf Definition.}\, Denote $\tau(a)$ by
$$\tau(a)=\tau(a,f,g)=\varliminf_{r\rightarrow\infty}\frac{\overline{N}_{0}(r,1/(f-a))}{\overline{N}(r,1/(f-a))}$$
where $\overline{N}_{0}(r,1/(f-a))$ means the common zeros of $f-a$ and $g-a$ with the same multiplicities, and  $\overline{N}(r,1/(f-a))\not\equiv0$. Otherwise, $\tau(a)=1$.
\par
\vskip 2mm E. Muse first applied {\bf Definition 3} to investigated the {\bf Open Question}. He proved the following theorem.
\par
\vskip 2mm{\bf Theorem E.}\, Suppose that $f$ and $g$ are two non-constant meromorphic functions that share four
values $a_{1}, a_{2}, a_{3}, a_{4}$. If $a_{1}$ is shared CM, then the conclusion that $f$ and $g$ share $a_{1}, a_{2}, a_{3}, a_{4}$ CM is true, provided that one of the following conditions is assumed in addition:\\
(i) $\tau(a_{2})\geq\frac{2}{3}$;\\
(ii) $\tau(a_{2})\geq\frac{1}{2}$ and the cross ratio $(a_{1},a_{2},a_{3},a_{4})=-1$.
\par
\vskip 2mm After the work of E. Muse  \cite{m1}, there were so many contributions on the uniqueness of meromorphic functions sharing four values. We list them in the following.
\par
\vskip 2mm {\bf Theorem E1.} \cite{wj,hb}\, Suppose that $f$ and $g$ are two non-constant meromorphic functions that share four values $a_{1}, a_{2}, a_{3}, a_{4}$. If $a_{1}$ is shared CM and $\tau(a_{2})\geq\frac{1}{2}$, then $f$ and $g$ share $a_{1}, a_{2}, a_{3}, a_{4}$ CM.
\par
\vskip 2mm{\bf Theorem E2.} \cite{re,ws,hb}\, Suppose that $f$ and $g$ are two non-constant meromorphic functions that share four
values $a_{1}, a_{2}, a_{3}, a_{4}$. Then the conclusion that $f$ and $g$ share $a_{1}, a_{2}, a_{3}, a_{4}$ CM is true, provided that one of the following conditions is assumed in addition:\\
(i) Two of the values satisfy $\tau(a_{i})\geq\frac{4}{5}$ for $i\in\{1,2,3,4\}$;\\
(ii) Two of the values satisfy $\tau(a_{i})\geq\frac{2}{3}$ for $i\in\{1,2,3,4\}$ and the cross ratio $(a_{1},a_{2},a_{3},a_{4})=-1$.
\par
\vskip 2mm {\bf Theorem E3.} \cite{re,sc}\, Suppose that $f$ and $g$ are two non-constant meromorphic functions that share four
values $a_{1}, a_{2}, a_{3}, a_{4}$. If three of the values satisfy $\tau(a_{i})\geq\frac{3}{4}$ for $i=1,2,3,4$, then $f$ and $g$ share $a_{1}, a_{2}, a_{3}, a_{4}$ CM.
\par
\vskip 2mm{\bf Theorem E4.} \cite{g3}\, Let $f$ and $g$ be two non-constant meromorphic functions, and let $a_{1}, a_{2}, a_{3}, a_{4}$ be four distinct constants. If $f$ and $g$ share $a_{1}$, $a_{2}$, $a_{3}$ IM, and $ a_{4}$ CM, and if  there exists a constant $\lambda>\frac{4}{5}$ and some set $I\subset(0,\infty)$ that has infinite linear measure such that
$$\frac{N(r,a_{4},f)}{T(r,f)}\geq \lambda$$
for all $r\in I$, then $f$ and $g$ share $a_{1}, a_{2}, a_{3}, a_{4}$ CM.
\par
\vskip 2mm{\bf Theorem E5.} \cite{re,u,wj}, Suppose that $f$ and $g$ are two non-constant meromorphic functions that share four values $a_{1}, a_{2}, a_{3}, a_{4}$. If all of the values satisfy $\tau(a_{i})\geq\frac{2}{3}$ for $i\in\{1,2,3,4\}$, then $f$ and $g$ share $a_{1}, a_{2}, a_{3}, a_{4}$ CM.
\par
\vskip 2mm{\bf Remark 1.1.}\, It is easy to check that the conclusion in  Theorem E4 is still true when $f$ and $g$ share $a_{1}$, $a_{2}$, $a_{3}$ IM* and $ a_{4}$ CM*.
\par
\vskip 2mm It does not seem to be a quite satisfactory answer to the "3IM+1CM" question,
for some restrictions. Based on {\bf Open Question}, it is interesting to study the uniqueness of meromorphic functions that share four small functions. So  according to  Theorem A1,  Theorem B1,  Theorem C1 and  Theorem D1, one may ask whether it is still true if we replace four values by four small functions.
\par
\vskip 2mm{\bf Conjecture 1.2.}\, If two non-constant meromorphic functions share three small functions IM and share a fourth small function CM, then do the functions necessarily share all four small functions CM?
\par
\vskip 2mm Obviously, {\bf Conjecture 1.2} implies {\bf Conjecture 1.2}.

\section{Meromorphic functions that share  small functions with their kth order derivatives}
\par
\vskip 2mm In 1977, Rubel-Yang \cite{ruy} first considered the uniqueness  about an entire function and its derivative. They proved:
\par
\vskip 2mm{\bf Theorem 3.1.} \ Let $f$ be a non-constant entire function, and let $a, b$ be two  distinct finite complex values. If $f$ and $f'$ share $a$ and $b$ CM, then $f\equiv f'$.
\par
\vskip 2mm And then,  Mues-Steinmetz \cite{muess}, Gundersen \cite{g0} improved Theorem 3.1   by the following result:
\par
\vskip 2mm{\bf Theorem 3.2.} \ Let $f$ be a non-constant entire function, and let $a, b$ be two  distinct finite complex values. If $f$ and $f'$ share $a$ and $b$ IM, then $f\equiv f'$.
\par
\vskip 2mm In 2000, Li-Yang  \cite{liy2} improved  Theorem 2.3. by proving following result:
\par
\vskip 2mm{\bf Theorem 3.3.}
Let $f$ be a non-constant entire function, let $k$ be a positive integer, and let $a$ and $b$ be distinct two finite values. If $f$ and $f^{(k)}$ share  $a$ and $b$ IM , then $f\equiv f^{(k)}$.
\par
\vskip 2mm In 2000, G. D. Qiu\cite{q} proved:
\par
\vskip 2mm {\bf Theorem 3.4.} Let $f$  be a  non-constant entire function,  and let $a\not\equiv\infty, b\not\equiv\infty$ be two  distinct small functions of $f$. If $f$ and $f'$ share $a$  and  $b$ IM, then  $f\equiv f'$.
\par
\vskip 2mm Li-Yang raised a conjecture in \cite{liy2}:
\par
\vskip 2mm {\bf Conjecture 3.1.}\quad Let $f(z)$ be a non-constant entire function, let $k$ be a positive integer, and let $a\not\equiv\infty, b\not\equiv\infty$ be two distinct small meromorphic functions  of $f$. If $f$ and $f^{(k)}$ share $a$ and $b$ IM, is $f\equiv f^{(k)}$?
\par
\vskip 2mm In 1992, Frank-Schwick \cite{fw} improved Theorem C and proved the following theorem:
\par
\vskip 2mm {\bf Theorem 3.5.}
Let $f$ be a non-constant meromorphic  function, let $k$ be a positive integer, and let $a, b, c$  be three distinct  values. If $f$ and $f^{(k)}$ share  $a$, $b$ and $c$ IM, then $f\equiv f^{(k)}$.
\par
\vskip 2mm It is natural to raise a conjecture:
\par
\vskip 2mm {\bf Conjecture 3.2}
Let $f$ be a non-constant meromorphic  function, let $k$ be a positive integer, and let $a, b, c$  be three distinct  small functions. If $f$ and $f^{(k)}$ share  $a$, $b$ and $c$ IM, is $f\equiv f^{(k)}$?

\par
\vskip 2mm {\bf Gol'dberg Conjecture.}\, Let $f$ be a transcendental meromorphic function in the complex plane and $k$ a positive integer, then
\[
\overline{N}(r,f)\leq N\left(r,\frac{1}{f^{(k)}}\right)+o(T(r,f)).
\]
To investigate above conjecture, we need a help of a deep result. In 2013, K. Yamanoi\cite{ya} proved the famous Gol’dberg conjecture, and in his paper, he obtained a more general result:
\par
\vskip 2mm {\bf Theorem 3.6.}\cite{ya}\, Let $f$ be a transcendental meromorphic function in the complex plane, let
$k\geq 1$ be an integer, let $\varepsilon >0$ and let $A $ be a finite set of finite complex numbers. Then we have
$$(k-1)\overline{N}(r,f)+\sum_{a\in A}N_{1}\left(r,\frac{1}{f-a}\right)\leq N\left(r,\frac{1}{f^{(k)}}\right)+\varepsilon T(r,f)$$
for all $r>e$ outside a set $E\subset(e,\infty)$  of logarithmic density $0$. Here $E$ depends on $f,A,k$ and $\varepsilon$, and the the notation of $N_{1}\left(r,\frac{1}{f-a}\right)$ is defined as in \eqref{eq4}.

\section{Uniqueness of Meromorphic and their shift or  n-exact order difference operators}

\par
\vskip 2mm  In recent years, there are many people work on the value distribution of meromorphic functions with respect to difference analogue,  see [2-5, 7-16, 18-22,27]. Heittokangas et al \cite{hkl} proved a similar result analogue of Theorem 3.1 concerning shift.
\par
\vskip 2mm {\bf Theorem 4.2.}\,
 Let $f(z)$ be a non-constant entire function of finite order, let $c$ be a nonzero finite complex value, and let $a, b$ be two finite distinct complex values.
If $f(z)$ and $f(z+\eta)$ share $a, b$ CM, then $f(z)\equiv f(z+c).$
\par
\vskip 2mm  Recently,  Chen-Yi  \cite{cy}, Zhang-Liao \cite{zl},  and Liu-Yang-Fang \cite{lyf}  proved:
\par
\vskip 2mm {\bf Theorem 4.3.}\,Let $f$ be a transcendental entire function of finite order,  let $c$ be a non-zero complex number, $n$ be a positive integer,  and let $a, b$ be two distinct small functions of $f$. If $ f$ and $\Delta_{c}^{n}f$ share $a$, $b$ CM, then $ f\equiv \Delta_{c}^{n}f$.
\par
\vskip 2mm  In 2021, the author \cite{hxh1,hxh2} considered the uniqueness of entire function that sharing one value CM and another value IM with its higher difference operators. He proved the following result.
\par
\vskip 2mm  {\bf Theorem 4.4.}\, Let $f(z)$ be a  transcendental entire function of finite order, let $c\neq0$ be a finite complex number, $n\geq1, k\geq0$  two  integers and let $a, b$ be two  distinct finite complex values. If $f(z)$ and $(\Delta_{c}^{n}f(z))^{(k)}$ share $a$ CM and share $b$ IM, then $a=2b$, and either $f(z)\equiv(\Delta_{c}^{n}f(z))^{(k)}$ or
 $$f(z)=be^{2(c_{1}z+d)}-2be^{c_{1}z+d}+2b,$$
$$(\Delta_{c}^{n}f(z))^{(k)}=be^{c_{1}z+d},$$
where $c_{1}=(-2)^{-\frac{n+1}{k}}$ for $k\geq1$ and $d$ are two finite constants.
\par
\vskip 2mm {\bf Remark 4.1}  Theorem D still holds when $f$ is a transcendental meromorphic function of finite order with $\overline{N}(r,f)=o(T(r,f))$.
\par
\vskip 2mm  In 2014, Halburd-Korhonen-Tohge \cite{h3} investigated the relationship of characteristic functions between $f(z)$ and $f(z+c)$ in $\rho_{2}(f)<1$. They obtain the following Lemma 5.12. Immediately, Theorem 4.1-Theorem 4.4 are still true when finite order is replaced by $\rho_{2}(f)<1$.
\par
\vskip 2mm  Li-Yi-Kang \cite{lyk}, L$\ddot{u}$-L$\ddot{u}$ \cite{ll}, Gao, et.al \cite{gkzz} improved Theorem C to meromorphic function. They proved.
\par
\vskip 2mm {\bf Theorem 4.5.}\, Let $f$ be a  transcendental meromorphic function of hyper-order less than $1$, let $c$ be a non-zero complex number, $n\geq1$, an integer, and let $a,b,c$  distinct values. If $f$ and $\Delta_{c}^{n}f$ share $a,b,c$ CM, then $f\equiv\Delta_{c}^{n}f$.
\par
\vskip 2mm In 2020, Chen-Xu \cite{cx} considered the uniqueness of meromorphic function that sharing $0,\infty$ CM and $1$ IM with its first order difference operator. They obtained the following theorem.
\par
\vskip 2mm {\bf Theorem 4.6.}\, Let $f$ be a  transcendental meromorphic function of hyper-order less than $1$, let $c$ be a non-zero complex number. If $f$ and $\Delta_{c}f$ share $0,\infty$ CM and $1$ IM, then $f\equiv\Delta_{c}f$.
\par
\vskip 2mm It is natural to raise a question from Theorem 4.6 that
\par
\vskip 2mm {\bf Conjecture 4.1.} Under the same conditions as in {\bf Theorem F}, can we improve first order difference operator to n-th order exact difference operators? And can we replace "$0,\infty$ CM and $1$ IM" by "$a,\infty$ CM and $b$ IM"? Here, $a$ and $b$ are two distinct finite values.
\par
\vskip 2mm In 2021, Huang-Fang \cite{hf} proved the following result.
\par
\vskip 2mm{\bf Theorem 4.7.}\, Let $f$ be a transcendental entire function with $\rho_{2}(f)<1$, let $c$ be a nonzero finite complex value, and
let $a, b$ be two distinct finite values. If $f'$ and $f_{c}$ share $a, b$ IM, then $f'\equiv f_{c}.$
\par
\vskip 2mm It is natural to ask a conjecture that:
\par
\vskip 2mm{\bf Conjecture 4.2.}\, As in Theorem 4.7, can we replace two distinct finite values by two distinct small functions?
\par
\vskip 2mm {\bf Conjecture 4.3.} \, Let $f$ be a transcendental entire function with $\rho_{2}(f)<1$, let $c$ be a nonzero finite complex value, $n$ be a positive integer and
let $a, b$ be two distinct small functions of $f$. If $f$ and $\Delta_{c}^{n}f$ share  $a,b$ IM, is $f\equiv\Delta_{c}^{n}f$?
\par
\vskip 2mm
{\bf Remark 4.2.}\quad It is interesting to investigate an analogue of shift version of  Theorem D, that is we consider the uniqueness of the case when $f^{(k)}$ and $f_{c}$  $a,b,d,\infty$ four distinct values IM. Here, we state the negative answer, i.e., $f^{(k)}\not\equiv f_{c}$. We explain the reason. $f_{c}$ and $f^{(k)}$ share poles IM (ignoring multiplicities). Suppose that $z_{0}$ is a pole of $f$, hence also $f^{(k)}(z_{0})=\infty$. Thus $f_{c}(z_{0})=f(nc+z_{0})=\infty$. By induction we obtain $f(c+z_{0})=\infty$ for all integers $n \geq 0$. Next, set $z_{1}=z_{0}-c$. Then $f_{c}(z_{1})=f(c+z_{1})=f(z_{0})=\infty$. Hence
$f^{(k)}(z_{1})=\infty$, so $f(z_{1})=\infty$, that is, $f(z_{0}-c)=\infty$. By induction we
see that $f(nc+z_{0})=\infty$ for all integers $n \leq0$ and hence for all integers
$n$.
\par
\vskip 2mm This means that the poles of $f$ must occur in sets of the form $nc+z_{0}$
where $n$ runs through all integers. This is, of course, the case if $f$ has
period $c$, but if $f$ does not have period $c$, it imposes a serious restriction
on what $f$ can look like even before the question of sharing other values
or functions is considered. In contrast, if we consider $f$ and $f^{(k)}$, then
the poles are automatically shared IM, so such a restriction does not
appear.
\par
\vskip 2mm Consider then if $f^{(k)}\equiv f_{c}$. Suppose that $f$ has at least one pole, say $z_{0}$, of order $m \geq 1$. Then $f^{(k)}$ has a pole of order $m + k$ at $z_{0}$. Define $z_{1}=z_{0}-c$. Since $f^{(k)}\equiv f_{c}$ and $f_{c}(z_{1})=f(c+z_{1})=f(z_{0})=\infty$, it follows that $f_{c}$ has a pole at $z_{1}$, so that $f^{(k)}$ has a  pole at $z_{1}$, and hence $f$ has a pole at $z_{1}$. Now $f_{c}$ has a pole of order $m$ at $z_{1}$, so $f^{(k)}$ has a pole of order $m$ at $z_{1}$, hence $f$ has a pole of order $m-k$ at $z_{1}$. If $m-k < 1$, we have a contradiction, and if $m-k \geq 1$, we may continue in this fashion ($f$ has a pole of order$m-2k$ at $z_{0}-2c$, and so on) and reach a contradiction after finitely many steps. Hence non-entire functions $f$ such that $f^{(k)}\equiv f_{c}$ do not exist at all.

\section{Meromorphic function that share pairs of small functions}
\par
\vskip 2mm Brosch \cite{bo}, Czubiak and Gundersen \cite{cg}, and Gundersen \cite{g}, Steinmetz \cite{sn}, Gundersen-Steinmetz-Tohge \cite{gst} studied shared pairs of values. Reinders \cite{r1} pointed out that $f$ and $g$ share the pairs of $(-\frac{1}{2},\frac{1}{4})$ CM, and thus $f$ and $g$ share five pairs of values, but $f$ is not a M$\ddot{o}$bius transformation of $g$. For five shared pairs of values, Gundersen etc. \cite{gst} proved the next result, which is the best possible:
\par
\vskip 2mm{\bf Theorem 5.1} \ Let $f$ and $g$ be two non-constant meromorphic functions,  and let $(a_{1},b_{1}), (a_{2},b_{2}), (a_{3},b_{3}), (a_{4},b_{4}), (a_{5},b_{5})$ be five pairs if values such that $a_{i}\neq a_{j}, b_{i}\neq b_{j}$ for $i,j\in\{1,2,3,4,5\}$ whenever $i\neq j$. If $f$ 
and $g$ share $(a_{1},b_{1}), (a_{2},b_{2})$ CM and share $(a_{3},b_{3}), (a_{4},b_{4},(a_{5}, b_{5})$ IM, then $f$ is a M$\rm \ddot{o}$bius transformation of $g$.
\par
\vskip 2mm In \cite{sn}, the author proved the following result.
\par
\vskip 2mm{\bf Theorem 5.2} \ Let $f$ and $g$ be two non-constant meromorphic functions,  and let $(a_{1},b_{1}), (a_{2},b_{2}), (a_{3},b_{3}), (a_{4},b_{4}), (a_{5},b_{5})$ be five pairs of values such that $a_{i}\neq a_{j}, b_{i}\neq b_{j}$ for $i,j\in\{1,2,3,4,5\}$ whenever $i\neq j$. If $f$ and $g$ share four pairs $(a_{1},b_{1}), (a_{2},b_{2}), (a_{3},b_{3}), (a_{4},b_{4}),$ and  $(a_{5}, b_{5})$ CM such that
$$m\left(r,\frac{1}{f-a_{5}}\right)+m\left(r,\frac{1}{g-b_{5}}\right)=S(r,f),$$
holds. Then either $f$ and $g$ are M$\rm \ddot{o}$bius transformations of each other or else $f=M_{1}\circ\widehat{f}\circ h$ and $g=M_{2}\circ\widehat{g}\circ h$
hold for suitably chosen M$\rm \ddot{o}$bius transformations $M_{1}$ and $M_{2}$ and some non-constant entire function $h$.
\par
\vskip 2mm Li and Qiao \cite{lq} proved that Theorem A is still valid for five distinct small functions, and they proved the following:
\par
\vskip 2mm {\bf Theorem 5.3} \ Let $f$ and $g$ be two non-constant meromorphic functions,
and let $a_{1}, a_{2}, a_{3}, a_{4}, a_{5}$  (one of them can be $\infty $) be five distinct small
functions of $f$ and $g$. If $f$ and $g$ share $a_{1}, a_{2}, a_{3}, a_{4}, a_{5}$ IM, then $f\equiv g$.
\par
\vskip 2mm Li and Yang \cite{liy3}, Zhang and Yang \cite{zy}, Nguyen and Si \cite{ns} studied shared pairs of small functions.  In 2009, Li and and Yang \cite{liy3} obtained three shared pairs of small functions CM and another shared pair of small functions IM theorem.  
\par
\vskip 2mm{\bf Theorem 5.4} Let $f$ and $g$ be two non-constant meromorphic functions, and let $(a_{1},b_{1}), (a_{2},b_{2}), (a_{3},b_{3}), (a_{4},b_{4})$ be four pairs of  small functions with respect to $f$ and $g$ such that $a_{i}\neq a_{j}, b_{i}\neq b_{j}$ for $i,j\in\{1,2,3,4\}$ whenever $i\neq j$. If $f$ and $g$ share  $(a_{1},b_{1}), (a_{2},b_{2}), (a_{3},b_{3})$ CM*, and $(a_{4},b_{4})$ IM*, then $f$ is a quasi-M$\ddot{o}$bius transformation of $g$.
\par
\vskip 2mm The following example showed that for four shared pairs of values, $2CM+2IM\neq4CM$.
\par
\vskip 2mm{\bf Example 5.1.} \cite{liy3}\quad Let $f=\frac{-(e^{z}-1)}{e^{z}-2}$ and $g=\frac{-2(e^{z}-1)^{2}}{e^{z}-2}$. By computation, we can verify that $f$ and $g$ share four values $(1,1), (\infty,\infty)$ CM, and $(0,0), (-2.-8)$, but $f$ is not a M$\ddot{o}$bius transformation of $g$.

\par
\vskip 2mm{\bf Conjecture 5.1.} \cite{liy3}\quad Let $f$ and $g$ be two non-constant meromorphic functions,  and let $(a_{1},b_{1}), (a_{2},b_{2}), (a_{3},b_{3}), (a_{4},b_{4}), (a_{5},b_{5})$ be five pairs of small functions with respect to $f$ and $g$ such that $a_{i}\neq a_{j}, b_{i}\neq b_{j}$ for $i,j\in\{1,2,3,4,5\}$ whenever $i\neq j$. If $f$ and $g$ share  $(a_{1},b_{1}), (a_{2},b_{2})$ CM and $(a_{3},b_{3}), (a_{4},b_{4}), (a_{5},b_{5})$ IM,  must $f$  be a quasi-M$\ddot{o}$bius transformation of $g$?

\section{Meromorphic functions that four small functions with their derivatives}
In 1990,  Yang L \cite{yl} first considered the derivatives of meromorphic functions sharing four finite values problem, and he proved the following theorem.
\par
\vskip 2mm {\bf Theorem 6.1} Let $f$ and $g$ be two non-constant meromorphic functions, and let $a_{i}(i=1,2,3,4)$ be four distinct finite constants. If $f^{(k)}$ and $g^{(k)}$ share $a_{i}(i=1,2,3,4)$ CM, where $k\geq1$ is a positive integer, then $f^{(k)}\equiv g^{(k)}$.
\par
\vskip 2mm In 2002, G.D. Qiu \cite{qgd2} proved the following result, which improved Theorem E
\par
\vskip 2mm {\bf Theorem 6.2} Let $f$ and $g$ be two non-constant meromorphic functions, and let $a_{i}\not\equiv\infty(i=1,2,3,4)$ be four distinct small functions of $f$ and $g$. If $f^{(k)}$ and $g^{(k)}$ share $a_{i}(i=1,2,3,4)$ CM, where $k\geq1$ is a positive integer, then   $f^{(k)}\equiv g^{(k)}$.
\par
\vskip 2mm In 2009, Li-Yi \cite{lyh1} improved Theorem F from four CM to four IM, and they proved the following result.
\par
\vskip 2mm{\bf Theorem 6.3} Let $f$ and $g$ be two non-constant meromorphic functions, and let $a_{i}\not\equiv\infty(i=1,2,3,4)$ be four distinct small functions of $f$ and $g$. If $f^{(k)}$ and $g^{(k)}$ share $a_{i}(i=1,2,3,4)$ IM, where $k\geq1$ is a positive integer,  then   $f^{(k)}\equiv g^{(k)}$.
\par
\vskip 2mm In 2011, Li-Yi-Hu \cite{lyh2} considered one of the four small functions of Theorem 6.3 is $\infty$. They proved the following theorem.
\par
\vskip 2mm {\bf Theorem 6.4} Let $f$ and $g$ be two non-constant meromorphic functions, and let $a_{1}, a_{2}, a_{3}$ be three distinct finite values. If $f^{(k)}$ and $g^{(k)}$ share $a_{1}, a_{2}, a_{3}, \infty$ IM, where $k\geq1$ is a positive integer, then $f$ and $g$ are functions with normal growth, have the same order of growth, a positive integer or infinity.   
\par
\vskip 2mm The authors raised a conjecture in the paper.
\par
\vskip 2mm{\bf Conjecture 6.1} Is the conclusion of Theorem 6.4 still true if we replace three distinct finite values by three distinct small functions?

\section{Main Results}
In this note, we obtain the following results.
\par
\vskip 2mm {\bf Theorem 1.}\, Let $f$ and $g$ be two non-constant meromorphic functions, and let $a_{1}, a_{2}, a_{3}, a_{4}$ be four distinct small functions of $f$ and $g$. If $f$ and $g$ share $a_{1}, a_{2}, a_{3}$ IM and $ a_{4}$ CM, then $a_{1}, a_{2}, a_{3}, a_{4}$ are four constants, and $f$ and $g$ share four distinct values CM.
\par
\vskip 2mm {\bf Theorem 2.}\, Let $f$ be a non-constant meromorphic  function, let $k$ be a positive integer, and let $a, b, c$  be three distinct  small functions such that none of them is identical to $\infty$. If $f$ and $f^{(k)}$ share  $a$, $b$ and $c$ IM, then $f\equiv f^{(k)}$.
\par
\vskip 2mm{\bf Theorem 3.} Let $f$ be a non-constant entire function, let $k$ be a positive integer, and let $a\not\equiv\infty, b\not\equiv\infty$ be two distinct small meromorphic functions  of $f$. If $f$ and $f^{(k)}$ share $a$ and $b$ IM, then $f\equiv f^{(k)}$.
\par
\vskip 2mm {\bf Theorem 4.}\, Let $f$ be a   meromorphic function of $\rho_{2}(f)<1$,  let $c\neq0$ be a finite complex number, $n\geq1, k\geq0$  two  integers and let $a, b$ be two distinct finite complex value. If $f$ and $(\Delta_{c}^{n}f)^{(k)}$ share $a,\infty$ CM and share $b$ IM, then either $f(z)\equiv(\Delta_{c}^{n}f(z))^{(k)}$ or
 $$f(z)=He^{2p}-ae^{p}+a,$$
$$(\Delta_{c}^{n}f(z))^{(k)}=He^{p},$$
where H is a meromorphic function with $\overline{N}(r,H)=o(T(r,f))$, and $p$ is a non-constant entire function of order less than $1$.

\

{\bf Theorem 5.}\, Let $f$ be a non-constant entire function of $\rho_{2}(f)<1$, $c\neq0$ a finite constant, let $n$ be a positive integer, and let  $a\not\equiv\infty$ and $b\not\equiv\infty$ be two distinct small functions of $f(z)$, if $f$ and $\Delta_{c}^{n}f$ share $a$ and $b$ IM, then $f\equiv \Delta_{c}^{n}f$.

\par
\vskip 2mm{\bf Theorem 6.}\, Let $f$ be a non-constant entire function of $\rho_{2}(f)<1$, let $c$ be a nonzero finite value, $k$ a positive integer, and let  $a\not\equiv\infty$ and $b\not\equiv\infty$ be two distinct small functions of $f$. If $f^{(k)}$ and $f_{c}$ share $a,b$ IM, then $f^{(k)}\equiv f_{c}$.
\par
\vskip 2mm {\bf Theorem 7.}\,  Let $f$ and $g$ be two non-constant meromorphic functions,  and let $(a_{1},b_{1}), (a_{2},b_{2}), (a_{3},b_{3}), (a_{4},b_{4}), (a_{5},b_{5})$ be five pairs of small functions with respect to $f$ and $g$ such that $a_{i}\neq a_{j}, b_{i}\neq b_{j}$ for $i,j\in\{1,2,3,4,5\}$ whenever $i\neq j$. If $f$ and $g$ share  $(a_{1},b_{1}), (a_{2},b_{2})$ CM and $(a_{3},b_{3}), (a_{4},b_{4}), (a_{5},b_{5})$ IM, then $f$ is a quasi-M$\ddot{o}$bius transformation of $g$.

\par
\vskip 2mm {\bf Theorem 8} let $f$ and $g$ be two non-constant meromorphic functions, and let $a_{i}(i=1,2,3)$ be three distinct small functions of $f$ and $g$. If $f^{(k)}$ and $g^{(k)}$ share three small functions $a_{1},a_{2},a_{3},\infty$ IM*, then $f^{(k)}$ and $g^{(k)}$ share four distinct small functions CM, and $f$ and $g$ are functions with normal growth, have the same order of growth, a positive integer or infinity.

\section{ Lemmas}
{\bf  Lemma 1} \cite{y1} Let $f$ be a non-constant meromorphic function, and let $a_{1}, a_{2}, a_{3}$ be three distinct small functions of $f$. Then
$$T(r,f)\leq \sum_{i=1}^{3}\overline{N}(r,\frac{1}{f-a_{i}})+o(T(r,f)).$$

{\bf  Lemma 2} \cite{y1}
Let $f$ and $g$ be two non-constant polynomials, and let $a, b$ be two distinct finite values. If $f$ and $g$ share $a, b$ IM, then $f\equiv g$.

{\bf  Lemma 3} \cite{ya} Let $f$ be a transcendental  meromorphic function, and let $v:R_{>e}\rightarrow\mathbb{N}_{>0} $ be a function such that
\begin{align}
v(r)\sim(\log^{+}\frac{T(r,f)}{\log r})^{20}.\label{2.1}
\end{align}
 Then we have
$$\bar{m}_{0,v(r)}(r,f)+\sum_{a\in\widehat{\mathbb{C}}}\hat{N}_{1}(r,a,f)=2T(r,f)+o(T(r,f)),$$
for all $r\to\infty$ outside a set $E$  of logarithmic density $0$.

\par
\vskip 2mm{\bf Lemma 4}  \cite{y3,y5,y6} Let $f$ be a non-constant meromorphic function, and let $a_{1}, a_{2}, a_{3},\ldots,a_{q}$ be $q$ distinct small functions of $f$. Then 
\begin{align}
(q-2-\varepsilon)T(r,f)\leq \sum_{i=1}^{q}\overline{N}\left(r,\frac{1}{f-a_{i}}\right)+O(1)\label{eq5}
\end{align}
for all $\varepsilon>0$ outside a set $E\subset(2,\infty)$ with $\int_{E}dr<\infty$. Especially, when all of small functions $a_{1}, a_{2}, a_{3},\ldots,a_{q}$ are $q$ distinct rational function, and when $f$ is a transcendental meromorphic function, the term $\varepsilon T(r,f)$ can be replace by $S(r,f)$, see \cite{y5}.

\par
\vskip 2mm{\bf Remark 6.1.} In \cite{y3}, the measure of the set $E$ in the first version of the truncated $q$-small function theorem  is $\int_{E}d\log\log r<\infty$, while in Page 2 of \cite{y4}, Yamanoi gave a correction to the measure of the set $E$, he estimated $E$ by  $\int_{E}\rm dr<\infty$.

{\bf  Lemma 5} Let $f$ and $g$ be two non-constant  meromorphic functions and let $a_{i}(i=1,2,3,4)$  be four distinct small functions of $f$ and $g$.  If $f$ and $g$ share $a_{i}(i=1,2,3,4)$ IM*, and if $f\not\equiv g$, then for some $\varepsilon>0$ and any small function $c\not\equiv a_{1}, a_{2}, a_{3}, a_{4}$ we have\\
(i):\quad  $T(r,f)= T(r,g)+\varepsilon T(r,f)$,\quad $T(r,f)= T(r,g)+\varepsilon T(r,f)$;\\
(ii):\quad $2T(r,f)=\sum_{i=1}^{4}\overline{N}(r,\frac{1}{f-a_{i}})+\varepsilon T(r,f)$;\\
(iii):\quad $\lim_{r\to\infty}\frac{\overline{N}(r,\frac{1}{f-c})}{T(r,f)}\geq1$.
\begin{proof}
\par
\vskip 2mm{\bf Case 1.}\quad None of $a_{j}$ for $j\in\{1,2,3,4\}$ is $\infty$. from the assumption that $f$ and $g$ share $a_{1}, a_{2}, a_{3}, a_{4}$ IM* and Lemma 2.1 we have
\begin{align}
&2T(r,f)\leq \sum_{j=1}^{4}\overline{N}\left(r,\frac{1}{f-a_{j}}\right)+\varepsilon T(r,f)\leq N\left(r,\frac{1}{f-g}\right)+\varepsilon T(r,g)\notag\\
&\leq T(r,f-g)+\varepsilon T(r,f)\leq T(r,f)+T(r,g)+\varepsilon T(r,f),\label{eq5}
\end{align}
as $r\to\infty$, possibly outside a set  $ E_{1}\subset(2,\infty)$  with $\int_{E_{1}}\rm dr<\infty$. From \eqref{eq5} we have
\begin{align}
T(r,f)\leq T(r,g)+\varepsilon T(r,f).\label{eq6}
\end{align}
Similarly,
\begin{align}
T(r,g)\leq T(r,f)+\varepsilon T(r,g),\label{eq7}
\end{align}
as $r\to\infty$, possibly outside a set  $ E_{2}\subset(2,\infty)$  with $\int_{E_{2}}\rm d r<\infty$. Therefore, from \eqref{eq6} and \eqref{eq7} we have $o(T(r,f))=o(T(r,g))$, as $r\not\in E_{1}\cup E_{2}$ and $r\to\infty$. From  \eqref{eq6} and \eqref{eq7} we have
\begin{align}
T(r,f)= T(r,g)+\varepsilon T(r,f).\label{eq8}
\end{align}
From \eqref{eq5} and \eqref{eq8} we have
\begin{align}
2T(r,f)&=\sum_{j=1}^{3}\overline{N}\left(r,\frac{1}{f-a_{j}}\right)+\overline{N}(r,f)+\varepsilon T(r,f),\label{eq9}
\end{align}
as $r\not\in E_{1}\cup E_{2}$ and $r\to\infty$.
\par
\vskip 2mm {\bf Case 2.}\quad One of $a_{j}$ for $j\in\{1,2,3,4\}$ is $\infty$, without loss of generality, we suppose $a_{4}=\infty$. First of all, we let $d$ be a finite complex constant such that $a_{j}\not\equiv d$ for $j\in\{1,2,3,4\}$. Next we set
$$F=\frac{1}{f-d},\quad G=\frac{1}{g-d},$$
then from the assumption that $f$ and $g$ share $a_{1}, a_{2},a_{3},\infty$ IM*, we obtain that $F$ and $G$ share $0$ and $d_{j}=\frac{1}{a_{j}-h}$ IM* for $j\in\{1,2,3\}$. with
the same method of {\bf Case 1}, we can obtain  \eqref{eq8} and 
$$2T(r,f)=\sum_{j=1}^{4}\overline{N}\left(r,\frac{1}{f-a_{j}}\right)+\varepsilon T(r,f),$$
as $r\to\infty$, possibly outside a set  $E_{1}\cup E_{2}\subset(2,\infty)$  with $\int_{E_{1}}\rm d r<\infty$ and $\int_{E_{2}}\rm d r<\infty$. Using Lemma 2.4 for any small function $c\not\equiv a_{1}, a_{2}, a_{3}, a_{4}$ we have
\begin{eqnarray*}
\begin{aligned}
3T(r,f)&\leq \overline{N}(r,f)+\overline{N}\left(r,\frac{1}{f-a_{1}}\right)+\overline{N}\left(r,\frac{1}{f-a_{2}}\right)+\overline{N}\left(r,\frac{1}{f-a_{3}}\right)\notag\\
&\quad+\overline{N}\left(r,\frac{1}{f-c}\right)+\varepsilon_{i} T(r,f)+o(T(r,f))\notag\\
&\leq (2+\varepsilon_{i})T(r,f)+\overline{N}\left(r,\frac{1}{f-c}\right)+o(T(r,f)),
\end{aligned}
\end{eqnarray*}
as $r\to\infty$, possibly outside a set  $ E_{1}\subset(2,\infty)$  with $\int_{E_{1}}\rm dr<\infty$. This deduces
 \begin{align}
T(r,f)\leq \overline{N}\left(r,\frac{1}{f-c}\right)+\varepsilon_{i} T(r,f)+o(T(r,f)),\label{eq11}
\end{align}
as $r\not\in E_{1}\cup E_{2}$ and $r\to\infty$. By \eqref{eq11} , we get
\begin{align}
\varliminf_{r\to\infty}\frac{\overline{N}\left(r,\frac{1}{f-c}\right)}{T(r,f)}\geq1,\label{eq12}
\end{align}
as $r\not\in E_{1}\cup E_{2}$ and $r\to\infty$.
\par
\vskip 2mm This proves Lemma 5.
\end{proof}

\par
\vskip 2mm {\bf  Lemma 6}\cite{hm} Suppose that $f$ is a transcendental meromorphic function and that $K$ is a positive number greater than one. Then there exists a set $M(K)$ of upper logarithmic density at most
$$\delta(K)=\min\{(2e^{K-1}-1)^{-1}, (1+e(K-1)\exp{e(1-K)})\}$$
such that the inequality
$$\limsup_{r\to\infty}\frac{T(r,f)}{T(r,f^{k})}\leq3eK,$$
holds for every positive integer $k$ and $r\not\in M(K)$.
\par
\vskip 2mm{\bf  Lemma 7}\cite{l} Let $F$ and $G$ be two distinct transcendental meromorphic functions such that $F$ and $G$ share $0, 1$ IM*, and let $b(\not\equiv0,1,\infty)$ be a small function of $F$ and $G$. If
$$\overline{N}(r,F)+\overline{N}(r,G)=\varepsilon(T(r,F)+T(r,G))$$
and
$$\overline{N}\left(r,\frac{1}{F-b}\right)+\overline{N}\left(r,\frac{1}{G-b}\right)-2\overline{N}_{E}(r,b)=\varepsilon(T(r,F)+T(r,G)),$$
where $\varepsilon$ is an arbitrary positive number, and $r\not\in E\cup E_{1}\cup M(K)$, of which $M(K)$ is defined in Lemma 6, then $F$ is a quasi-M$\ddot{o}$bius transformation of $G$.
\par
\vskip 2mm{\bf  Lemma 8}\cite{lyh1}  Let $f$ and $g$ be two non-constant meromorphic functions. If $f$ and $g$ share $0,1,\infty$ IM*, and $f$ is  a M$\ddot{o}$bius transformation of $g$,  then $f$ and $g$ assume one of the following six relations: (i) $fg=1$; (ii) $(f-1)(g-1)=1$; (iii) $f+g=1$; (iv) $f=cg$; (v) $f-1=c(g-1)$; (vi) $[(c-1)f+1][(c-1)g-c]=-c$, where $c\not\equiv0,1,\infty$ is a small function of $f$ and $g$.
\par
\vskip 2mm In order to prove "3IM+1CM", we need some Lemmas.
\par
\vskip 2mm{\bf Lemma 9} \cite{y1}\, Let $f$ and $g$ be two non-constant rational functions. If $f$ and $g$ share four distinct values $a_{1}, a_{2}, a_{3}, a_{4}$ IM, then $f\equiv g$.
\par
\vskip 2mm{\bf Lemma 10} \cite{y1}\, Let $f$ and $g$ be two non-constant  meromorphic functions and let $a_{1}, a_{2}, a_{3}, a_{4}$  be four distinct values.  Suppose we have
$$\phi=\frac{f'g'(f-g)^{2}}{(f-a_{1})(f-a_{2})(f-a_{3})(g-a_{1})(g-a_{2})(g-a_{3})}.$$
If $f$ and $g$ share $a_{i}(i=1,2,3,4)$ IM, and if $f\not\equiv g$, then\\
(i):\quad  $T(r,f)=T(r,g)+S(r,f), T(r,g)=T(r,f)+S(r,g)$;\\
(ii):\quad $2T(r,f)=\sum_{i=1}^{4}\overline{N}\left(r,\frac{1}{f-a_{i}}\right)+S(r,f)$\\
(iii):\quad $T(r,f)=\overline{N}\left(r,\frac{1}{f-b}\right)+S(r,f)$ for any $b\neq a_{i}(i=1,2,3,4)$.\\
(iv):\quad $N_{0}\left(r,\frac{1}{f'}\right)=S(r,f)$ and $N_{0}\left(r,\frac{1}{g'}\right)=S(r,f)$, where $N_{0}\left(r,\frac{1}{f'}\right)$ means the zeros of $f'$ but not the zeros of $f-a_{i}$ for $i=1,2,3,4$.\\
(v)\quad $N_{(m,n)}\left(r,\frac{1}{f-a_{i}}\right)=S(r,f)$, where $m\geq2,n\geq2$ are two positive integers, and $N_{(m,n)}(r,\frac{1}{f-a_{i}})$ is the counting function which count the  common zeros of $f-a_{i}$ with multiplicities $m$ and $g-a_{i}$ with multiplicities $n$.
\par
\vskip 2mm 
{\bf Lemma 11}\cite{y1} Let $f$ and $g$ be two non-constant meromorphic functions. If $f$ and $g$ share $0,1,\infty$ CM and $f\not\equiv g$, then
\begin{align}
N_{(2}(r,f)+N_{(2}\left(r,\frac{1}{f}\right)+N_{(2}\left(r,\frac{1}{f-1}\right)=S(r),\label{eq2.1}
\end{align}
where $N_{(2}(r,f)$ and $N_{(2}\left(r,\frac{1}{f-a}\right)$ denote the counting functions of the 
multiple poles and multiple $a$-points of $f$, respectively, and $S(r):=S(r,f)=S(r,g)$.

\par
\vskip 2mm{\bf Remark 6.2.} We can see from the proof of Lemma 2.1 in \cite{y1}(Theorem 5.4, Page 267) is still true if $f$ and $g$ share $0,1,\infty$ CM*.

{\bf Lemma 12}\cite{h3} Let $f$ be a non-constant meromorphic function with $\rho_{2}(f)<1$,  and let $c$ be a non-zero complex number. Then
$$m(r,\frac{f_{c}}{f})=o(T(r,f)),$$
for all r outside of a possible exceptional set E with finite logarithmic measure.

{\bf Lemma 13}\label{23l} Let $f_{1}$ and $f_{2}$ be  non-constant meromorphic functions in $|z|<\infty$, then
$$N(r,f_{1}f_{2})-N(r,\frac{1}{f_{1}f_{2}})=N(r,f_{1})+N(r,f_{2})-N(r,\frac{1}{f_{1}})-N(r,\frac{1}{f_{2}}),$$
where $0<r<\infty$.

{\bf Lemma 14}\label{23l} Let $f$ be a transcendental entire function, let $k$ be a  positive integer, and let $a\not\equiv\infty, b\not\equiv\infty$ be two distinct small meromorphic functions of $f$. Suppose
\[L(f_{c})=\left|\begin{array}{rrrr}a-b& &f_{c}-a \\
a'-b'& &f_{c}'-a'\end{array}\right|\]
and
\[L(f^{(k)})=\left|\begin{array}{rrrr}a-b& &f^{(k)}-a \\
a'-b'& &f^{(k+1)}-a'\end{array}\right|,\]
and $f_{c}$ and $f^{(k)}$ share $a$  and  $b$ IM,  then $L(f_{c})\not\equiv0$ and $L(f^{(k)})\not\equiv0$.

\begin{proof}
Suppose that $L(f_{c})\equiv0$, then we can get $\frac{f'_{c}-a'}{f_{c}-a}\equiv\frac{a'-b'}{a-b}$. Integrating both side of above we can obtain $f_{c}-a=C_{1}(a-b)$, where $C_{1}$ is a nonzero constant. Then we have $T(r,f_{c})=o(T(r,f))$, a contradiction. Hence $L(f_{c})\not\equiv0$.

Since $f^{(k)}$ and $f_{c}$ share $a$ and $b$ IM, and that $f$ is a non-constant entire function, then by Lemma 12, we get
\begin{align}
T(r,f_{c})&\leq \overline{N}(r,\frac{1}{f_{c}-a})+\overline {N}(r,\frac{1}{f_{c}-b})+o(T(r,f))\notag\\
&= \overline {N}(r,\frac{1}{f^{(k)}-a})+\overline {N}(r,\frac{1}{f^{(k)}-b})+o(T(r,f))\notag\\
&\leq 2T(r,f^{(k)})+o(T(r,f)).
\end{align}
Hence $a$ and $b$ are small functions of $f^{(k)}$. If $L(f^{(k)})\equiv0$, then we can get $f^{(k)}-a=C_{2}(a-b)$, where $C_{2}$ is a nonzero constant. And we get $T(r,f^{(k)})=o(T(r,f))$. Combing (2.1) we obtain $T(r,f_{c})=o(T(r,f))$, a contradiction.
\end{proof}

{\bf Lemma 15}\label{24l}  Let $f$ be a transcendental entire function, and $k$ a positive integer. Let $a\not\equiv\infty, b\not\equiv\infty$ be two distinct small meromorphic functions of $f$.  Again let $d_{j}=a-j(a-b)$, where $j\neq0,1$ is a positive integer. Then
$$m(r,\frac{L(f_{c})}{f_{c}-a})=o(T(r,f)), \quad m(r,\frac{L(f_{c})}{f_{c}-b})=o(T(r,f)).$$
And
$$m(r,\frac{L(f_{c})f_{c}}{(f_{c}-a)(f_{c}-b)(f_{c}-d_{j})})=o(T(r,f)),$$
where $L(f)$ is defined as in Lemma 2.14.

\begin{proof}
Obviously, we have
$$m(r,\frac{L(f_{c})}{f_{c}-a})\leq m(r,-\frac{(a'-b')(f_{c}-a)}{f_{c}-a})+m(r,\frac{(a-b)(f'_{c}-a')}{f_{c}-a})=o(T(r,f)),$$
and
$$\frac{L(f_{c})f_{c}}{(f_{c}-a)(f_{c}-b)(f_{c}-d_{j})}=\frac{C_{1}L(f_{c})}{f_{c}-a}+\frac{C_{2}L(f_{c})}{f_{c}-b}+\frac{C_{3}L(f_{c})}{f_{c}-d_{j}},$$
where $C_{i}(i=1,2,3)$ are small functions of $f$. Thus  we have
$$m(r,\frac{L(f_{c})f_{c}}{(f_{c}-a)(f_{c}-b)(f_{c}-d_{j})})=o(T(r,f)).$$
\end{proof}

{\bf Lemma 16}\label{21l} Let $f$ be a nonconstant meromorphic function with $\rho_{2}(f)<1$,  and let $c$ be a non-zero complex number. Then
$$T(r,f)=T(r,f_{c})+o(T(r,f)).$$

{\bf Lemma 17}\cite{zy} Let $f$ and $g$ be two non-constant meromorphic functions,  and let $(a_{1},b_{1}), (a_{2},b_{2})$ and $ (a_{3},b_{3})  (a_{4},b_{4}), (a_{5},b_{5})$ be five pairs of small functions with respect to $f$ and $g$ such that $a_{i}\neq a_{j}, b_{i}\neq b_{j}$ for $i,j\in\{1,2,3,4,5\}$ whenever $i\neq j$. If $f$ and $g$ share  $(a_{1},b_{1}), (a_{2},b_{2}), (a_{3},b_{3}), (a_{4},b_{4}), (a_{5},b_{5})$  IM*, and if  $f$ is not a quasi-M$\ddot{o}$bius transformation of $g$, then, for any $\varepsilon>0$, the following inequalities hold for any three distinct indices $i, j, k$ among the integers $\{1, 2, 3, 4, 5\}$:\\
$$\overline{N}\left(r,\frac{1}{f-a_{i}}\right)+\overline{N}\left(r,\frac{1}{f-a_{j}}\right)\leq \frac{3}{2}T(r,f)+S^{*}(r,f);$$
$$\overline{N}\left(r,\frac{1}{f-a_{i}}\right)+\overline{N}\left(r,\frac{1}{f-a_{j}}\right)+\overline{N}\left(r,\frac{1}{f-a_{k}}\right)\leq 2T(r,f)+S^{*}(r,f);$$
$$\frac{3}{2}T(r,f)\leq\overline{N}\left(r,\frac{1}{f-a_{i}}\right)+\overline{N}\left(r,\frac{1}{f-a_{j}}\right)+\overline{N}\left(r,\frac{1}{f-a_{k}}\right)+S^{*}(r,f);$$
$$\overline{N}\left(r,\frac{1}{f-a_{i}}\right)\leq 2\overline{N}\left(r,\frac{1}{f-a_{j}}\right)+S(r,f).$$

{\bf Lemma 18} \cite{liy3}  Let $f$ and $g$ be two non-constant meromorphic functions,  and let $(a_{1},b_{1}), (a_{2},b_{2})$ and $ (a_{3},b_{3})  (a_{4},b_{4}), (a_{5},b_{5})$ be five pairs of small functions with respect to $f$ and $g$ such that $a_{i}\neq a_{j}, b_{i}\neq b_{j}$ for $i,j\in\{1,2,3,4,5\}$ whenever $i\neq j$. If $f$ and $g$ share  $(a_{1},b_{1}), (a_{2},b_{2}), (a_{3},b_{3}), (a_{4},b_{4}), (a_{5},b_{5})$  IM*, and if  $f$ is not a quasi-M$\ddot{o}$bius transformation of $g$, then for $i\neq j$, we have
$$T(r,f)=T(r,g)+S^{*}(r,f);$$
$$T(r,f)\leq \overline{N}\left(r,\frac{1}{f-a_{i}}\right)+\overline{N}\left(r,\frac{1}{f-a_{j}}\right)+S^{*}(r,f);$$
$$3T(r,f)=\sum_{i=1}^{5}\overline{N}\left(r,\frac{1}{f-a_{i}}\right)+S^{*}(r,f).$$$c\not\equiv0,1,\infty$ is a small function of $f$ and $g$.

\section{The proof of Theorem 1 }
\par
\vskip 2mm If $f\equiv g$, there is nothing we need to prove. So we assume that $f\not\equiv g$.  If $f$ and $g$ are two rational functions,  then invoking Lemma 9, we have  $f\equiv g$, a contradiction. Hence $f$ and $g$ are transcendental. In the following,  we just  consider  Theorem 1 excluding the case  at least two CM shared small functions.
\par
\vskip 2mm{\bf Case 1} All of $a_{1}, a_{2}, a_{3}, a_{4}$ are constants. Without loss of generality in assuming that $a_{1}=0, a_{2}=1, a_{3}=b, a_{4}=\infty$, and that $\infty$ is the only CM* shared value.

Suppose we have
\begin{align}
F=f-A,\quad G=g-A,\label{eq7.1}
\end{align}
where $A$ is a non-constant rational function. Since $f$ and $g$ share $0,1,b$ IM* and $\infty$ CM*, we know that $F$ and $G$ share $-A,1-A,b-A$ IM* and $\infty$ CM*. Let $v(r)\ll p\leq\infty$ be a large enough positive integer, let $j=1,2,3,\ldots,p$ be an integer,  and let  $d_{j}\not\equiv 0,1,b,\infty$ are distinct values.   Applying Lemma 4 to $F$ for every $j\in\{1,2,3,\ldots,p\}$, we can obtain 
\begin{eqnarray*}
\begin{aligned}
&(p+2)T(r,F)\leq \overline{N}\left(r,f\right)+\overline{N}\left(r,\frac{1}{F+A}\right)+\overline{N}\left(r,\frac{1}{F-1+A}\right)+\overline{N}\left(r,\frac{1}{F-b+A}\right)\\
&\quad+\sum_{j=1}^{p}\overline{N}\left(r,\frac{1}{F-d_{j}}\right)+S(r,F)\leq 4T(r,F)+\sum_{j=1}^{p}\overline{N}\left(r,\frac{1}{F-d_{j}}\right)\\
&-m(r,F)-m\left(r,\frac{1}{F+A}\right)-m\left(r,\frac{1}{F-1+A}\right)-m\left(r,\frac{1}{F-b+A}\right)+S(r,F),
\end{aligned}
\end{eqnarray*}
as $r\to\infty$, possibly outside a set  $ E_{1}\subset(2,\infty)$  with $\int_{E_{1}}\rm dr<\infty$. This deduces
\begin{align}
&(p-2)T(r,F)+m(r,F)+m\left(r,\frac{1}{F+A}\right)+m\left(r,\frac{1}{F-1+A}\right)\notag\\
&+m\left(r,\frac{1}{F-b+A}\right)\leq \sum_{j=1}^{p}\overline{N}\left(r,\frac{1}{F-d_{j}}\right)+S(r,F),\label{eq7.2}
\end{align}
as $r\to\infty$ and $r\not\in E_{1}$. 
\par
\vskip 2mm Let $q=v(r)$ be a positive integer,  and as $v(r)\sim(\log^{+}\frac{T(r,f)}{\log r})^{20}=o(T(r,f))$, we have $\{1,2,3,\ldots,q\}\subset\{1,2,3,\ldots,p\}$ and $\{d_{1},d_{2},\ldots,d_{q}\}\subset\{d_{1},d_{2},\ldots,d_{p}\}$. Suppose $\hat{\mathbb{C}}=\{d_{1},d_{2},\ldots,d_{p}\}\cup \{a_{3}, a_{4}\}$. According to $\{d_{1},d_{2},\ldots,d_{q}\}\subset\{d_{1},d_{2},\ldots,d_{p}\}$, \eqref{eq1}, \eqref{eq2}, \eqref{eq4}, \eqref{eq7.2}, (v) of Lemma 3 and Lemma 5, we have
\begin{eqnarray*}
\begin{aligned}
&2T(r,F)=\bar{m}_{0,q}(r,F)+\sum_{d\in\hat{\mathbb{C}}}\hat{N}_{1}\left(r,\frac{1}{F-d_{j}}\right)+S(r,F)\leq m\left(r,F\right)+N\left(r,F\right)\\
&+\sum_{j=1}^{q}m\left(r,\frac{1}{F-d_{j}}\right)
-\overline{N}\left(r,F\right)+\sum_{j=1}^{p}N\left(r,\frac{1}{F-d_{j}}\right)
-\sum_{j=1}^{p}\overline{N}\left(r,\frac{1}{F-d_{j}}\right)+S(r,F)\\
&=(p+1)T(r,F)+\sum_{j=1}^{q}m\left(r,\frac{1}{F-d_{j}}\right)-\sum_{j=1}^{p}m\left(r,\frac{1}{F-d_{j}}\right)-(p-2)T(r,F)\\
&-\overline{N}\left(r,f\right)-m(r,F)-m\left(r,\frac{1}{F+A}\right)-m\left(r,\frac{1}{F-1+A}\right)-m\left(r,\frac{1}{F-b+A}\right)+S(r,F)\\
&\leq 2T(r,F)-m\left(r,\frac{1}{F+A}\right)-m\left(r,\frac{1}{F-1+A}\right)-m\left(r,\frac{1}{F-b+A}\right)+S(r,F)
\end{aligned}
\end{eqnarray*}
as $r\not\in E_{1}\cup E_{2}$ and $r\to\infty,$ where and in what follows, $E_{2}\subset(0,\infty)$ is a set of logarithmic density $0$. This and \eqref{eq7.2} deduce that
\begin{align}
m\left(r,\frac{1}{f}\right)+m\left(r,\frac{1}{f-1}\right)+m\left(r,\frac{1}{f-b}\right)\leq S(r,f),\label{eq7.3}
\end{align}
as  $r\not\in E_{1}\cup E_{2}$ and $r\to\infty$. Similarly, we also have
\begin{align}
m\left(r,\frac{1}{g}\right)+m\left(r,\frac{1}{g-1}\right)+m\left(r,\frac{1}{g-b}\right)\leq S(r,f)\leq S(r,g)=S(r,f),\label{eq7.4}
\end{align}
as $r\not\in E_{3}\cup E_{4}$ and $r\to\infty,$ where and in what follows, $ E_{3}\subset(2,\infty)$  with $\int_{E_{3}}\rm dr<\infty$ and $E_{4}\subset(0,\infty)$ is a set of logarithmic density $0$.
\par
\vskip 2mm Let
\begin{align}
F_{1}=\frac{1}{f},\quad G_{1}=\frac{1}{g}.\label{eq7.5}
\end{align}
Because  $f$ and $g$ share $0,1,b$ IM* and $\infty$ CM*, we know that $F_{1}$ and $G_{1}$ share $c=\frac{1}{b},1, \infty$ IM*, and $0$ CM*. Suppose we have the following two functions,
\begin{align}
\varphi=\frac{F'_{1}(F_{1}-G_{1})}{F_{1}(F_{1}-1)(F_{1}-c)},\label{eq7.6}
\end{align}
\begin{align}
\psi=\frac{G'_{1}(F_{1}-G_{1})}{G_{1}(G_{1}-1)(G_{1}-c)}.\label{eq7.7}
\end{align}
It is easy to see that $\varphi\not\equiv0$ and $\psi\not\equiv0$ since $F_{1}\not\equiv G_{1}$. On the one hand,
\begin{align}
m(r,\varphi)&=m\left(r,\frac{F'_{1}(F_{1}-G_{1})}{F_{1}(F_{1}-1)(F_{1}-c)}\right)\notag\\
&\leq m\left(r,\frac{F'_{1}F_{1}}{F_{1}(F_{1}-1)(F_{1}-c)}\right)+m\left(r,\frac{F'_{1}}{F_{1}(F_{1}-1)(F_{1}-c)}\right)+m(r,G_{1})\notag\\
&\leq m\left(r,\frac{1}{g}\right)+S(r,f)=S(r,f),\label{eq7.8}
\end{align}
as $r\not\in E_{1}\cup E_{2}\cup E_{3}\cup E_{4}$ and $r\to\infty$. Similarly, we have
\begin{align}
m(r,\psi)\leq m\left(r,\frac{1}{f}\right)+S(r,f)=S(r,f),\label{eq7.9}
\end{align}
as $r\not\in E_{1}\cup E_{2}\cup E_{3}\cup E_{4}$ and $r\to\infty$. Next, we analyze the poles of $\varphi$ and $\psi$. As  $F_{1}$ and $G_{1}$ share $1,c,\infty$ IM* and $0$ CM*, so we can see from \eqref{eq7.6} and \eqref{eq7.7} that
\begin{align}
N(r,\varphi)\leq \sum_{n=2}N_{(1,n)}(r,G_{1})-\sum_{n=2}\overline{N}_{(1,n)}(r,G_{1}),\label{eq7.10}
\end{align}
\begin{align}
N(r,\psi)\leq \sum_{m=2}N_{(m,1)}(r,F_{1})-\sum_{m=2}\overline{N}_{(m,1)}(r,F_{1}),\label{eq7.11}
\end{align}
where $m$ and $n$ are two positive integers.
\par
\vskip 2mm Note that the poles of $F$ with multiplicities $m\geq2$ and the poles of $G$ with multiplicities $1$ are the zero points of $\varphi$, so we get that
\begin{align}
&\sum_{m=2}N_{(m,1)}(r,F_{1})\leq N\left(r,\frac{1}{\varphi}\right)\leq T(r,\varphi)=N(r,\varphi)+S(r,f)\notag\\
&\leq \sum_{n=2}N_{(1,n)}(r,G_{1})-\sum_{n=2}\overline{N}_{(1,n)}(r,G_{1})+S(r,f),\label{eq7.12}
\end{align}
as $r\not\in E_{1}\cup E_{2}\cup E_{3}\cup E_{4}$ and $r\to\infty$. Similarly,
\begin{align}
&\sum_{n=2}N_{(1,n)}(r,G_{1})\leq N\left(r,\frac{1}{\psi}\right)\leq T(r,\psi)=N(r,\psi)+S(r,f)\notag\\
&\leq\sum_{m=2}N_{(m,1)}(r,F_{1})-\sum_{m=2}\overline{N}_{(m,1)}(r,F_{1})+S(r,f),\label{eq7.13}
\end{align}
as $r\not\in E_{1}\cup E_{2}\cup E_{3}\cup E_{4}$ and $r\to\infty$. Then \eqref{eq7.12} and \eqref{eq7.13} deduce
\begin{align}
\overline{N}(r,F_{1})\leq N_{(1,1)}(r,F_{1})+S(r,f),\label{eq7.14}
\end{align}
as $r\not\in E_{1}\cup E_{2}\cup E_{3}\cup E_{4}$ and $r\to\infty$. Likely, we also have
\begin{align}
\overline{N}\left(r,\frac{1}{F_{1}-c}\right)\leq N_{(1,1)}\left(r,\frac{1}{F_{1}-c}\right)+S(r,f)\label{eq7.15}
\end{align}
and
\begin{align}
\overline{N}\left(r,\frac{1}{F_{1}-1}\right)&\leq N_{(1,1)}\left(r,\frac{1}{F_{1}-1}\right)+S(r,f),\label{eq7.16}
\end{align}
as $r\not\in E_{1}\cup E_{2}\cup E_{3}\cup E_{4}$ and $r\to\infty$. If $\varphi\equiv\psi$, then $F_{1}$ and $G_{1}$ share four values CM, and hence $f$ and $g$ share four values CM. Hence, we assume that $\varphi\not\equiv\psi$. So according to a simple computation, we can know by \eqref{eq7.6} and \eqref{eq7.7} that
\begin{align}
& N_{(1,1)}\left(r,\frac{1}{F_{1}}\right)+N_{(1,1)}\left(r,\frac{1}{F_{1}-c}\right)+N_{(1,1)}\left(r,\frac{1}{F_{1}-1}\right)\leq N\left(r,\frac{1}{\varphi-\psi}\right)\notag\\
&\leq T(r,\varphi-\psi)= m(r,\varphi-\psi)+N(r,\varphi-\psi)\leq N(r,\varphi-\psi)+S(r,f)\notag\\
&\leq \sum_{n=2}N_{(1,n)}(r,G_{1})-\sum_{n=2}\overline{N}_{(1,n)}(r,G_{1})+\sum_{m=2}N_{(m,1)}(r,F_{1})-\sum_{m=2}\overline{N}_{(m,1)}(r,F_{1})\notag\\
&+S(r,f)=N(r,F_{1})-\overline{N}(r,F_{1})+S(r,f),\label{eq7.17}
\end{align}
as $r\not\in E_{1}\cup E_{2}\cup E_{3}\cup E_{4}$ and $r\to\infty$. Thus,  we apply (ii) of Lemma 5 and \eqref{eq7.14}-\eqref{eq7.17} that
\begin{align}
&2T(r,F_{1})=\overline{N}(r,F_{1})+\overline{N}\left(r,\frac{1}{F_{1}}\right)+\overline{N}\left(r,\frac{1}{F_{1}-c}\right)+\overline{N}\left(r,\frac{1}{F_{1}-1}\right)+S(r,F)\notag\\
&\leq N_{(1,1)}\left(r,\frac{1}{F_{1}-1}\right)+N_{(1,1)}\left(r,\frac{1}{F_{1}-c}\right)+N_{(1,1)}\left(r,\frac{1}{F_{1}-d}\right)+\overline{N}(r,F_{1})\notag\\
&+S(r,f)\leq N(r,F_{1})+S(r,f)\leq T(r,F_{1})+S(r,f),\label{eq7.18}
\end{align}
as $r\not\in E_{1}\cup E_{2}\cup E_{3}\cup E_{4}$ and $r\to\infty$. Thus, we obtain from \eqref{eq7.18} that 
\begin{align}
T(r,F_{1})=S(r,f),\label{eq7.19}
\end{align}
as $r\not\in E_{1}\cup E_{2}\cup E_{3}\cup E_{4}$ and $r\to\infty$. This is a contradiction. Therefore,  we obtain that $F_{1}$ and $G_{1}$ share four values CM, and hence $f$ and $g$ share four values CM.
\par
\vskip 2mm{\bf Case 2}\quad  Two of $a_{i}$ are non-constant small functions  for $i=1,2,3,4$, says  $a_{1}$ and $a_{2}$,  since we can do some transforms on $f$ and $g$. For example, we set
$$F=\frac{f-a_{1}}{f-a_{2}}\frac{a_{3}-a_{2}}{a_{3}-a_{1}},\quad G=\frac{g-a_{1}}{g-a_{2}}\frac{a_{3}-a_{2}}{a_{3}-a_{1}}.$$
If $\frac{a_{4}-a_{2}}{a_{4}-a_{1}}\frac{a_{3}-a_{2}}{a_{3}-a_{1}}\neq0,1,\infty$ is a constant, we can get from {\bf Case 1} and  the assumption $f$ and $g$ sharing $a_{1}, a_{2}, a_{3}$ IM* and $a_{4}$ CM*  that  $F$ and $G$ share $0,1,\frac{a_{4}-a_{2}}{a_{4}-a_{1}}\frac{a_{3}-a_{2}}{a_{3}-a_{1}},\infty$ CM*, and hence, $f$ and $g$ share $a_{1}, a_{2}, a_{3}, a_{4}$ CM. 
\par
\vskip 2mm If $\frac{a_{4}-a_{2}}{a_{4}-a_{1}}\frac{a_{3}-a_{2}}{a_{3}-a_{1}}$ is a non-constant small function, we set 
$$F_{1}=F-\frac{a_{4}-a_{2}}{a_{4}-a_{1}}\frac{a_{3}-a_{2}}{a_{3}-a_{1}},\quad G_{1}=G-\frac{a_{4}-a_{2}}{a_{4}-a_{1}}\frac{a_{3}-a_{2}}{a_{3}-a_{1}}.$$
Easy to see from above that $F_{1}$ and $G_{1}$ share $-\frac{a_{4}-a_{2}}{a_{4}-a_{1}}\frac{a_{3}-a_{2}}{a_{3}-a_{1}},1-\frac{a_{4}-a_{2}}{a_{4}-a_{1}}\frac{a_{3}-a_{2}}{a_{3}-a_{1}}, \infty$ IM* and $0$ CM*. If $F_{1}$ is a quasi-M$\ddot{o}$bius transformation of $G_{1}$, then 
$F$ is a quasi-M$\ddot{o}$bius transformation of $G$. Thus, $f$ and $g$ share $a_{1}, a_{2}, a_{3}, a_{4}$ CM. So we suppose $F_{1}$ is not a quasi-M$\ddot{o}$bius transformation of $G_{1}$. Without loss of generality, we only need to consider the case that $a_{1}$ and $ a_{2}$ are two non-constant  small functions, and  $f$ and $g$ share $a_{1}, a_{2}, a_{3}$ IM* and $a_{4}$ CM*. Let $v(r)\ll p\leq\infty$ be a large enough positive integer, let $j=1,2,3,\ldots,p$ be an integer,  and let  $d_{j}\not\equiv a_{1}, a_{2}, a_{3}, a_{4}$ are distinct values.   By Lemma 4 and (ii) of Lemma 5, we can obtain 
\begin{eqnarray*}
\begin{aligned}
&3T(r,f)\leq \overline{N}\left(r,\frac{1}{f-a_{1}}\right)+\overline{N}\left(r,\frac{1}{f-a_{2}}\right)+\overline{N}\left(r,\frac{1}{f-a_{3}}\right)+\overline{N}\left(r,\frac{1}{f-a_{4}}\right)\\
&\quad+\overline{N}\left(r,\frac{1}{f-d_{j}}\right)+\varepsilon_{j} T(r,f)+S(r,f)\leq (2+\varepsilon_{j})T(r,f)+\overline{N}\left(r,\frac{1}{f-d_{j}}\right)+S(r,f),
\end{aligned}
\end{eqnarray*}
as $r\to\infty$, possibly outside a set  $ E_{1}\subset(2,\infty)$  with $\int_{E_{1}}\rm dr<\infty$. This deduces
\begin{align}
T(r,f)\leq \overline{N}\left(r,\frac{1}{f-d_{j}}\right)+\varepsilon_{j} T(r,f)+S(r,f),\label{eq7.20}
\end{align}
as $r\to\infty$ and $r\not\in E_{1}$. By \eqref{eq7.20} and that $\varepsilon_{j}>0$ is any positive number for every $j\in\{1,2,3,\ldots,p\}$, we get
\begin{align}
\lim_{r\to\infty}\frac{\overline{N}\left(r,\frac{1}{f-d_{j}}\right)}{T(r,f)}\geq1,\label{eq7.21}
\end{align}
as $r\not\in E_{1}$ and $r\to\infty$.
\par
\vskip 2mm Let $q=v(r)$ be a positive integer,  and as $v(r)\sim(\log^{+}\frac{T(r,f)}{\log r})^{20}=o(T(r,f))$, we have $\{1,2,3,\ldots,q\}\subset\{1,2,3,\ldots,p\}$ and $\{d_{1},d_{2},\ldots,d_{q}\}\subset\{d_{1},d_{2},\ldots,d_{p}\}$. Suppose $\hat{\mathbb{C}}=\{d_{1},d_{2},\ldots,d_{p}\}\cup \{a_{3}, a_{4}\}$. According to $\{d_{1},d_{2},\ldots,d_{q}\}\subset\{d_{1},d_{2},\ldots,d_{p}\}$, \eqref{eq1}-\eqref{eq4}, \eqref{eq7.20}, (i)-(ii) of Lemma 5 and Lemma 3, we have
\begin{eqnarray*}
\begin{aligned}
&2T(r,f)=\bar{m}_{0,q}(r,f)+\sum_{d\in\hat{\mathbb{C}}}\hat{N}_{1}\left(r,\frac{1}{f-d_{j}}\right)+S(r,f)\leq m\left(r,\frac{1}{f-a_{3}}\right)+m\left(r,\frac{1}{f-a_{4}}\right)\\
&\quad+\sum_{j=1}^{q}m\left(r,\frac{1}{f-d_{j}}\right)+N\left(r,\frac{1}{f-a_{3}}\right)
+N\left(r,\frac{1}{f-a_{4}}\right)-\overline{N}\left(r,\frac{1}{f-a_{3}}\right)\\
&\quad-\overline{N}\left(r,\frac{1}{f-a_{4}}\right)+\sum_{j=1}^{p}N\left(r,\frac{1}{f-d_{j}}\right)
-\sum_{j=1}^{p}\overline{N}\left(r,\frac{1}{f-d_{j}}\right)+S(r,f)\\
&=(p+2)T(r,f)+\sum_{j=1}^{q}m\left(r,\frac{1}{f-d_{j}}\right)-\sum_{j=1}^{p}m\left(r,\frac{1}{f-d_{j}}\right)-\sum_{j=1}^{p}\overline{N}\left(r,\frac{1}{f-d_{j}}\right)\\
&\quad-\overline{N}\left(r,\frac{1}{f-a_{3}}\right)-\overline{N}\left(r,\frac{1}{f-a_{4}}\right)+S(r,f)\\
&\leq (p+2)T(r,f)-\sum_{j=1}^{p}\overline{N}\left(r,\frac{1}{f-d_{j}}\right)-\overline{N}\left(r,\frac{1}{f-a_{3}}\right)-\overline{N}\left(r,\frac{1}{f-a_{4}}\right)+S(r,f)
\end{aligned}
\end{eqnarray*}
as $r\not\in E_{1}\cup E_{2}$ and $r\to\infty,$ where and in what follows, $E_{2}\subset(0,\infty)$ is a set of logarithmic density $0$. This implies that
\begin{align}
\overline{N}\left(r,\frac{1}{f-a_{3}}\right)+\overline{N}\left(r,\frac{1}{f-a_{4}}\right)+\sum_{j=1}^{p}\overline{N}\left(r,\frac{1}{f-d_{j}}\right)\leq p T(r,f)+S(r,f),\label{eq7.22}
\end{align}
as $r\not\in E_{1}\cup E_{2}$ and $r\to\infty$. Divide $T(r,f)$ and we take the limite of the both sides of \eqref{eq7.22}, and combining \eqref{eq7.20} we can get
 \begin{align}
&\lim_{r\to\infty}\frac{\overline{N}\left(r,\frac{1}{f-a_{3}}\right)+\overline{N}\left(r,\frac{1}{f-a_{4}}\right)}{T(r,f)}+p\notag\\
&\leq \lim_{r\to\infty}\frac{\overline{N}\left(r,\frac{1}{f-a_{3}}\right)+\overline{N}\left(r,\frac{1}{f-a_{4}}\right)}{T(r,f)}+\sum_{j=1}^{p}\lim_{r\to\infty}\frac{\overline{N}\left(r,\frac{1}{f-d_{j}}\right)}{T(r,f)}\leq p,\label{eq7.23}
\end{align}
  as $r\to\infty$  and  $r\not\in E_{1}\cup E_{2}$. From \eqref{eq7.23} we have
\begin{eqnarray*}
\begin{aligned}
\lim_{r\to\infty}\frac{\overline{N}\left(r,\frac{1}{f-a_{3}}\right)+\overline{N}\left(r,\frac{1}{f-a_{4}}\right)}{T(r,f)}\leq0,
\end{aligned}
\end{eqnarray*}
as $r\to\infty$  and  $r\not\in E_{1}\cup E_{2}$. That is to say
\begin{align}
\overline{N}\left(r,\frac{1}{f-a_{3}}\right)+\overline{N}\left(r,\frac{1}{f-a_{4}}\right)\leq 2\varepsilon T(r,f),\label{eq7.24}
\end{align}
as $r\not\in E_{1}\cup E_{2}$ and $r\to\infty$. We set
$$F_{2}=\frac{f-a_{1}}{f-a_{3}}\frac{a_{2}-a_{3}}{a_{2}-a_{1}}, \quad G_{2}=\frac{g-a_{1}}{g-a_{3}}\frac{a_{2}-a_{3}}{a_{2}-a_{1}}.$$
As $f$ and $g$ share $a_{1}, a_{2}, a_{3}$ IM*, and $a_{4}$ CM*, we know that $F_{2}$ and $G_{2}$ share $0,1,\infty$ IM* and $c=\frac{a_{2}-a_{3}}{a_{2}-a_{1}}\frac{a_{4}-a_{1}}{a_{4}-a_{3}}$ CM*. If $c$ is a constant, then by {\bf Case 1}, we deduce that $F_{2}$ and $G_{2}$ share $0,1,c,\infty$ CM, and hence $f$ and $g$ share $a_{1}, a_{2}, a_{3}, a_{4}$ CM. If $c\neq0,1,\infty$ is a non-constant small function, then   \eqref{eq7.1} can deduce
\begin{align}
\overline{N}(r,F_{2})+\overline{N}\left(r,\frac{1}{F_{2}-c}\right)+\overline{N}(r,G_{2})+\overline{N}\left(r,\frac{1}{G_{2}-c}\right)\leq 2\varepsilon T(r,f),\label{eq7.25}
\end{align}
as $r\not\in E_{1}\cup E_{2}\cup E_{4}$ and $r\to\infty$. Therefore, combining \eqref{eq7.25} and Lemma 7, we can get $F_{2}$ is a quasi-M$\ddot{o}$bius transformation of $G_{2}$. Then by Lemma 8, we can yield that $F_{2}$ and $G_{2}$ share $0,1,\infty$ CM* and $c$ IM*. Thus,  by Theorem D1, we obtain $F_{2}$ and $G_{2}$ share four values CM. That is to say, $f$ and $g$ share $a_{1}, a_{2}, a_{3}, a_{4}$ CM.

This completes Theorem 1.

\section{Proof of Theorem 2}
\par
\vskip 2mm The proof of Theorem 2 can be found in \cite{hxh3}.

\section{The proof of Theorem 3 }
\par
\vskip 2mm We prove Theorem 3 by contradiction. Suppose that $f\not\equiv f^{(k)}$. If $f$ is a non-constant polynomial, then $a$ and $b$ must be two finite values. By Lemma 2, we obtain $f\equiv f^{(k)}$, a contradiction. So $f$ is transcendental, since $f$ is a transcendental entire function, and the assumption that $f$ and $f^{(k)}$ share $a$ and $b$ IM, from the First Fundamental Theorem and Lemma 1  we have
\begin{eqnarray*}
\begin{aligned}
T(r,f)&\leq \overline{N}(r,\frac{1}{f-a})+\overline{N}(r,\frac{1}{f-b})+o(T(r,f))\\
&= \overline{N}(r,\frac{1}{f^{(k)}-a})+\overline{N}(r,\frac{1}{f^{(k)}-b})+o(T(r,f))\\
&\leq N(r,\frac{1}{f-f^{(k)}})+o(T(r,f))\\
&\leq T(r,f-f^{(k)})+o(T(r,f))= m(r,f-f^{(k)})+o(T(r,f))\\
&\leq m(r,f)+m(r,1-\frac{f^{(k)}}{f})+o(T(r,f))\\
&\leq T(r,f)+o(T(r,f)),
\end{aligned}
\end{eqnarray*}
as $r\to\infty$, possibly outside a set  $ E_{1}\subset(2,\infty)$  with $\int_{E_{1}}\rm dr<\infty$. It follows that
\begin{eqnarray}
T(r,f)=\overline{N}(r,\frac{1}{f-a})+\overline{N}(r,\frac{1}{f-b})+o(T(r,f)),\label{eq9.1}
\end{eqnarray}
as  $r\not\in E_{1}$ and $r\to\infty$.
\par
\vskip 2mm We discuss the following three cases.
\par
\vskip 2mm {\bf Case 1}\quad $a$ and $b$ are finite constants. By Theorem 3.4, we have $f\equiv f^{(k)}$, a contradiction. 
\par
\vskip 2mm{\bf Case 2}\quad $a$ and $b$ are two non-constant small functions. Let $p\leq\infty$ be a large enough positive integer, let $i=1,2,3,\ldots,p$ be an integer, and let $h_{i}$ be $p$ distinct finite values. From Lemma 4 and \eqref{eq9.1} we have
\begin{eqnarray*}
\begin{aligned}
&2T(r,f)\leq \overline{N}(r,f)+\overline{N}\left(r,\frac{1}{f-a}\right)+\overline{N}\left(r,\frac{1}{f-b}\right)\notag\\
&+\overline{N}\left(r,\frac{1}{f-h_{i}}\right)+\varepsilon_{i}T(r,f)\leq T(r,f)+\overline{N}\left(r,\frac{1}{f-h_{i}}\right)+\varepsilon_{i}T(r,f)+o(T(r,f)),
\end{aligned}
\end{eqnarray*}
as $r\to\infty$, possibly outside a set  $ E_{1}\cup E_{3}\subset(2,\infty)$  with $\int_{E_{3}}\rm dr<\infty$. This deduces that 
\begin{align}
\lim_{r\to\infty}\frac{\overline{N}\left(r,\frac{1}{f-h_{i}}\right)}{T(r,f)}\geq1,\label{eq9.2}
\end{align}
 as $r\not\in E_{1}\cup E_{2}$ and $r\to\infty$, and where $\varepsilon_{i}>0$ is any positive number for every $i\in\{1,2,3,\ldots,p\}$.  Let $q=v(r)$ be a positive integer,  since $v(r)\sim\left(\log^{+}\frac{T(r,f)}{\log r}\right)^{20}=o(T(r,f))$, we have $\{1,2,3,\ldots,q\}\subset\{1,2,3,\ldots,p\}$ and $\{h_{1},h_{2},\ldots,h_{q}\}\subset\{h_{1},h_{2},\ldots,h_{p}\}$.  Suppose $\hat{\mathbb{C}}=\{h_{1},h_{2},\ldots,h_{j}\}\cup \{\infty\}$. According to $\{h_{1},h_{2},\ldots,h_{q}\}\subset\{h_{1},h_{2},\ldots,h_{p}\}$, \eqref{eq1}, \eqref{eq2}, \eqref{eq4} and Lemma 2.3
\begin{eqnarray*}
\begin{aligned}
&2T(r,f)=\bar{m}_{0,q}(r,f)+\sum_{h\in\hat{\mathbb{C}}}\hat{N}_{1}\left(r,\frac{1}{f-h_{i}}\right)+o(T(r,f))\\
&\leq m(r,f)+\sum_{i=1}^{q}m\left(r,\frac{1}{f-h_{i}}\right)+\sum_{i=1}^{p}N\left(r,\frac{1}{f-h_{i}}\right)-\sum_{i=1}^{p}\overline{N}\left(r,\frac{1}{f-h_{i}}\right)+o(T(r,f))\\
&\leq (p+1)T(r,f)+\sum_{i=1}^{q}m\left(r,\frac{1}{f-h_{i}}\right)-\sum_{i=1}^{p}m\left(r,\frac{1}{f-h_{i}}\right)-\sum_{i=1}^{p}\overline{N}\left(r,\frac{1}{f-h_{i}}\right)\\
&\quad+o(T(r,f))\leq (p+1)T(r,f)-\sum_{i=1}^{p}\overline{N}\left(r,\frac{1}{f-h_{i}}\right)+\varepsilon T(r,f),
\end{aligned}
\end{eqnarray*}
as $r\not\in E_{1}\cup E_{3}\cup E_{4}$ and $r\to\infty,$ where and in what follows, $E_{4}\subset(0,\infty)$ is a set of logarithmic density $0$. This yields
\begin{align}
T(r,f)+\sum_{i=1}^{p}\overline{N}\left(r,\frac{1}{f-h_{i}}\right)\leq pT(r,f)+o(T(r,f)),\label{eq9.3}
\end{align}
 as $r\to\infty$  and  $r\not\in E_{1}\cup E_{3}\cup E_{4}$.   Divide $T(r,f)$ and we take the limit of the both sides of \eqref{eq9.3}, and combining \eqref{eq9.2} we can get
\begin{eqnarray*}
\begin{aligned}
1+p\leq 1+\sum_{i=1}^{p}\lim_{r\to\infty}\frac{\overline{N}\left(r,\frac{1}{f-h_{i}}\right)}{T(r,f)}\leq p,
\end{aligned}
\end{eqnarray*}
as $r\to\infty$  and  $r\not\in E_{1}\cup E_{3}\cup E_{4}$. This implies that $1\leq0$, a contradiction.
\par
\vskip 2mm{\bf Case 3}\quad $a$ is a non-constant small function, and $b$ is a finite value. Let 
\begin{align}
 F=\frac{f-a}{b-a}-c,\quad G=\frac{f^{(k)}-a}{b-a}-c, \label{eq9.4}
\end{align}
where $c\not\equiv 0,1$ is a non-constant small function. As $f$ and $f^{(k)}$ share $a$ and $b$ IM, we can know from \eqref{eq3.4} that $F$ and $G$ share two non-constant small functions $-c$ and $1-c$ IM*. From \eqref{eq9.4} we have
\begin{align}
&\overline{N}(r,\frac{1}{F+c})=\overline{N}(r,\frac{1}{f-a})+o(T(r,f))\notag\\
&\overline{N}(r,\frac{1}{F-1+c})=\overline{N}(r,\frac{1}{f-b})+o(T(r,f))\notag\\
&\overline{N}(r,F)+\overline{N}(r,f)+o(T(r,f)).\label{eq9.5}
\end{align}

and we also obtain from \eqref{eq9.1} and \eqref{eq9.5} that
\begin{align}
T(r,F)&=T(r,f)+o(T(r,f))=\overline{N}(r,\frac{1}{f-a})+\overline{N}(r,\frac{1}{f-b})+o(T(r,f))\notag\\
&=\overline{N}(r,\frac{1}{F+c})+\overline{N}(r,\frac{1}{F-1+c})+o(T(r,f)).\label{eq9.6}
\end{align}
Then similar to the proof of {\bf Case 2}, we can obtain a contradiction.
\par
\vskip 2mm This completes Theorem 3.

\section{The proof of Theorem 4}
\par
\vskip 2mm We prove by contradiction. Assume that $f\not\equiv (\Delta_{c}^{n}f(z))^{(k)}$. Since $f$ is a meromorphic function of $\rho_{2}(f)<1$,  $f$ and $(\Delta_{\eta}^{n}f)^{(k)}$ share $a, \infty$ CM, then  we get
\begin{align}
\frac{(\Delta_{\eta}^{n}f)^{(k)}-a}{f-a}=e^{p},\label{eq3.1}
\end{align}
where $p$ is an entire function of order less than $1$. If $f$ is a non-constant rational function, we obtain from \eqref{eq3.1} that $e^{p}$ is a non-zero finite constant. If $e^{p}=1$, we have $f\not\equiv (\Delta_{c}^{n}f(z))^{(k)}$, a contradiction. Thus, $e^{p}\neq1$, and from the assumption that $f$ and $(\Delta_{\eta}^{n}f)^{(k)}$ share $b$ IM, we obtain
\begin{align}
N\left(r,\frac{1}{f-b}\right)=N\left(r,\frac{1}{(\Delta_{\eta}^{n}f)^{(k)}-b}\right)=0.\label{eq3.2}
\end{align}
That is to say, $f$ and $(\Delta_{\eta}^{n}f)^{(k)}$ share $a, b, \infty$ CM. Hence, there exists a constant $h\neq1$ such that
\begin{align}
\frac{(\Delta_{\eta}^{n}f)^{(k)}-b}{f-b}=e^{h}.\label{eq3.3}
\end{align}
Similarly,
\begin{align}
N\left(r,\frac{1}{f-a}\right)=N\left(r,\frac{1}{(\Delta_{\eta}^{n}f)^{(k)}-a}\right)=0.\label{eq3.4}
\end{align}
On the other hand, by \eqref{eq3.1} and Lemma 12, we get
\begin{align}
m\left(r,\frac{1}{f-a}\right)=m(r,e^{p})=o(T(r,f)),\label{eq3.5}
\end{align}
for all r outside of a possible exceptional set $E_{1}$ with finite logarithmic measure. From \eqref{eq3.4} and \eqref{eq3.5} we have $T(r,f)=o(T(r,f))$, a contradiction. Therefore, $f$ is a transcendental meromorphic function of $\rho_{2}(f)<1$.
\par
\vskip 2mm If $\overline{N}(r,f)=o(T(r,f))$, from {\bf Remark 4.1} and the same method of proving Theorem 4.4 we have either $f(z)\equiv(\Delta_{c}^{n}f(z))^{(k)}$, a contradiction. Or
 $$f(z)=He^{2p}-ae^{p}+a,$$
$$(\Delta_{c}^{n}f(z))^{(k)}=He^{p},$$
where H is a meromorphic function with $\overline{N}(r,H)=o(T(r,f))$, and $p$ is a non-constant entire function of order less than $1$. So, we assume $\overline{N}(r,f)\neq o(T(r,f))$ in the following.
\par
\vskip 2mm Suppose
\begin{align}
F=f-A,\quad   G=(\Delta_{c}^{n}f(z))^{(k)}-A,  \label{eq9.6}
\end{align}
where $A\not\equiv a,b,\infty$ is a non-constant rational function. From the fact that $f$ and $(\Delta_{\eta}^{n}f)^{(k)}$ share $a, \infty$ CM and $b$ IM, we obtain from \eqref{eq9.6} that $F$ and $G$ share $a-A, \infty$ CM* and $b-A$ IM*. Let $j\leq\infty$ be a large enough positive integer, let $i=1,2,3,\ldots,j$ be an integer, and let $d_{i}\not\equiv a-A,b-A,\infty$ be $j$ distinct finite values.  By  Lemma 4, we get
\begin{eqnarray*}
\begin{aligned}
\quad &(j+1)T(r,F)\\
&\leq \overline{N}(r,F)+\overline{N}\left(r,\frac{1}{F-a+A}\right)+\overline{N}\left(r,\frac{1}{F-b+A}\right)+\sum_{i=1}^{j}\overline{N}\left(r,\frac{1}{F-d_{i}}\right)+o(T(r,F))\\
&\leq 3T(r,F)-m(r,F)-m\left(r,\frac{1}{F-a+A}\right)-m\left(r,\frac{1}{F-b+A}\right)+\sum_{i=1}^{j}\overline{N}\left(r,\frac{1}{F-d_{i}}\right)\\
&\quad+o(T(r,F))\leq 3T(r,F)-m(r,F)\\
&\quad-m\left(r,\frac{1}{F-a+A}\right)-m\left(r,\frac{1}{F-b+A}\right)+\sum_{i=1}^{j}\overline{N}\left(r,\frac{1}{F-d_{i}}\right)+o(T(r,F)),
\end{aligned}
\end{eqnarray*}
as $r\to\infty$, possibly outside a set  $ E_{1}\cup E_{2}$  with $\int_{E_{2}}\rm dr<\infty$. This deduces
\begin{align}
&\quad (j-2)T(r,F)+m(r,F)+m\left(r,\frac{1}{F-a+A}\right)+m\left(r,\frac{1}{F-b+A}\right)\notag\\
&\leq \sum_{i=1}^{j}\overline{N}\left(r,\frac{1}{F-d_{i}}\right)+o(T(r,F)),\label{eq3.7}
\end{align}
 as $r\not\in E_{1}\cup E_{2}$ and $r\to\infty$. Let $q=v(r)$ be a positive integer,  let $l\in\{1,2,3,\ldots,q\}$ be an integer, and let $d_{l}$ be $l$ distinct finite values. As $v(r)\sim\left(\log^{+}\frac{T(r,F)}{\log r}\right)^{20}=o(T(r,F))$, we have $\{1,2,3,\ldots,q\}\subset\{1,2,3,\ldots,j\}$ and $\{d_{1},d_{2},\ldots,d_{q}\}\subset\{d_{1},d_{2},\ldots,d_{j}\}$. Suppose $\hat{\mathbb{C}}=\{d_{1},d_{2},\ldots,d_{j}\}\cup \{\infty\}$. According to $\{d_{1},d_{2},\ldots,d_{q}\}\subset\{d_{1},d_{2},\ldots,d_{j}\}$, \eqref{eq1}, \eqref{eq2}, \eqref{eq4}, \eqref{eq3.7},  and Lemma 3, we have
\begin{eqnarray*}
\begin{aligned}
&2T(r,F)=\bar{m}_{0,q}(r,F)+\sum_{d_{i}\in\hat{\mathbb{C}}}\hat{N}_{1}\left(r,\frac{1}{F-d_{i}}\right)+o(T(r,F))\\
&\leq \sum_{l=1}^{q}m\left(r,\frac{1}{F-d_{l}}\right)+m(r,F)+N(r,F)-\overline{N}(r,F)+\sum_{i=1}^{j}N\left(r,\frac{1}{F-d_{i}}\right)\\
&\quad-\sum_{i=1}^{j}\overline{N}\left(r,\frac{1}{F-d_{i}}\right)+o(T(r,F))=T(r,F)-\overline{N}(r,F)\\
&\quad+\sum_{l=1}^{q}m\left(r,\frac{1}{F-d_{i}}\right)+\sum_{i=1}^{j}N\left(r,\frac{1}{F-d_{i}}\right)-(j-2)T(r,F)\\
&\quad-m(r,F)-m\left(r,\frac{1}{F-a+A}\right)-m\left(r,\frac{1}{F-b+A}\right)+o(T(r,F))\\
&\leq 3T(r,F)-\overline{N}(r,F)-m(r,F)-m\left(r,\frac{1}{F-a+A}\right)-m\left(r,\frac{1}{F-b+A}\right)\\
&+\sum_{l=1}^{q}m\left(r,\frac{1}{F-d_{l}}\right)-\sum_{i=1}^{j}m\left(r,\frac{1}{F-d_{i}}\right)+o(T(r,F))\\
&\leq 3T(r,F)-\overline{N}(r,F)-m(r,F)-m\left(r,\frac{1}{F-a+A}\right)-m\left(r,\frac{1}{F-b+A}\right)+o(T(r,F)),
\end{aligned}
\end{eqnarray*}
as $r\not\in E_{1}\cup E_{2}\cup E_{3}$ and $r\to\infty,$ where and in what follows, $E_{3}\subset(0,\infty)$ is a set of logarithmic density $0$. This yields that
\begin{align}
&\overline{N}(r,F)+m(r,F)+m\left(r,\frac{1}{F-a+A}\right)+m\left(r,\frac{1}{F-b+A}\right)\notag\\
&\leq T(r,F)+\varepsilon T(r,F),\label{eq9.8}
\end{align}
  as $r\to\infty$  and  $r\not\in E_{1}\cup E_{2}\cup E_{3}$. By \eqref{eq3.1}, \eqref{eq9.6} and Lemma 12, we get
\begin{align}
m\left(r,\frac{1}{F-a+A}\right)&=m\left(r,\frac{1}{f-a}\right)+o(T(r,f))\notag\\
&=m(r,e^{p})+o(T(r,f)),\label{eq3.9}
\end{align}
as $r\not\in E_{1}$ and $r\to\infty$. Rewrite \eqref{eq3.1} as
\begin{align}
\frac{(\Delta_{\eta}^{n}f)^{(k)}-f}{f-a}=e^{p}-1.\label{eq3.10}
\end{align}
From \eqref{eq3.10} we have
\begin{align}
\overline{N}\left(r,\frac{1}{f-b}\right)\leq\overline{N}\left(r,\frac{1}{e^{p}-1}\right).\label{eq3.11}
\end{align}
Applying the second fundamental theorem of Nevanlinna to $e^{p}$, we get
\begin{eqnarray*}
\begin{aligned}
&m(r,e^{p})=T(r,e^{p})\leq\overline{N}(r,e^{p})+\overline{N}\left(r,\frac{1}{e^{p}}\right)+\overline{N}\left(r,\frac{1}{e^{p}-1}\right)+o(T(r,e^{p}))\notag\\
&\leq \overline{N}\left(r,\frac{1}{e^{p}-1}\right)+o(T(r,e^{p}))\leq T(r,e^{p})+o(T(r,e^{p})),
\end{aligned}
\end{eqnarray*}
as $r\not\in E_{1}$ and $r\to\infty$. This implies that
\begin{align}
m(r,e^{p})=T(r,e^{p})=\overline{N}\left(r,\frac{1}{e^{p}-1}\right)+o(T(r,e^{p})).\label{eq3.12}
\end{align}
From \eqref{eq3.6}, \eqref{eq9.8}-\eqref{eq3.9} and \eqref{eq3.11}-\eqref{eq3.12} we have
\begin{align}
\overline{N}(r,F)+\overline{N}\left(r,\frac{1}{F-b+A}\right)\leq T(r,F)+o(T(r,F)),\label{eq3.13}
\end{align}
as $r\not\in E_{1}\cup E_{2}\cup E_{3}$ and $r\to\infty$. By  Lemma 4 and \eqref{eq3.13}, we get
\begin{eqnarray*}
\begin{aligned}
\quad &(j+1)T(r,F)\\
&\leq \overline{N}(r,F)+\overline{N}\left(r,\frac{1}{F-a+A}\right)+\overline{N}\left(r,\frac{1}{F-b+A}\right)+\sum_{i=1}^{j}\overline{N}\left(r,\frac{1}{F-d_{i}}\right)+o(T(r,F))\\
&\leq 2T(r,F)-m\left(r,\frac{1}{F-a+A}\right)+\sum_{i=1}^{j}\overline{N}\left(r,\frac{1}{F-d_{i}}\right)\\
&\quad+o(T(r,F))\leq 2T(r,F)-m\left(r,\frac{1}{F-a+A}\right)+\sum_{i=1}^{j}\overline{N}\left(r,\frac{1}{F-d_{i}}\right)+o(T(r,F)),
\end{aligned}
\end{eqnarray*}
 as $r\not\in E_{1}\cup E_{2}\cup E_{3}$ and $r\to\infty$. This deduces that
\begin{align}
&\quad (j-1)T(r,F)+m\left(r,\frac{1}{F-a+A}\right)\notag\\
&\leq \sum_{i=1}^{j}\overline{N}\left(r,\frac{1}{F-d_{i}}\right)+o(T(r,F)),\label{eq3.14}
\end{align}
 as $r\not\in E_{1}\cup E_{2}\cup E_{3}$ and $r\to\infty$. According to $\{d_{1},d_{2},\ldots,d_{q}\}\subset\{d_{1},d_{2},\ldots,d_{j}\}$, \eqref{eq1}, \eqref{eq2}, \eqref{eq4}, \eqref{eq3.14},  and Lemma 3, we have
\begin{eqnarray*}
\begin{aligned}
&2T(r,F)=\bar{m}_{0,q}(r,F)+\sum_{d_{i}\in\hat{\mathbb{C}}}\hat{N}_{1}\left(r,\frac{1}{F-d_{i}}\right)+o(T(r,F))\\
&\leq \sum_{l=1}^{q}m\left(r,\frac{1}{F-d_{l}}\right)+m(r,F)+N(r,F)-\overline{N}(r,F)+\sum_{i=1}^{j}N\left(r,\frac{1}{F-d_{i}}\right)\\
&\quad-\sum_{i=1}^{j}\overline{N}\left(r,\frac{1}{F-d_{i}}\right)+o(T(r,F))=T(r,F)-\overline{N}(r,F)\\
&\quad+\sum_{l=1}^{q}m\left(r,\frac{1}{F-d_{i}}\right)+\sum_{i=1}^{j}N\left(r,\frac{1}{F-d_{i}}\right)-(j-1)T(r,F)\\
&\quad-m\left(r,\frac{1}{F-a+A}\right)+2\varepsilon T(r,F)\leq T(r,F)-\overline{N}(r,F)-m\left(r,\frac{1}{F-a+A}\right)\\
&\quad+\sum_{l=1}^{q}m\left(r,\frac{1}{F-d_{l}}\right)-\sum_{i=1}^{j}m\left(r,\frac{1}{F-d_{i}}\right)\\
&\leq 2T(r,F)-\overline{N}(r,F)-m\left(r,\frac{1}{F-a+A}\right)+o(T(r,F)),
\end{aligned}
\end{eqnarray*}
as $r\not\in E_{1}\cup E_{2}\cup E_{3}\subset(0,\infty)$ and $r\to\infty$. This with \eqref{eq9.6}, \eqref{eq3.9} and \eqref{eq3.12} deduce that
\begin{align}
\overline{N}(r,F)+\overline{N}\left(r,\frac{1}{F-b+A}\right)\leq o(T(r,F)),\label{eq3.15}
\end{align}
as $r\not\in E_{1}\cup E_{2}\cup E_{3}$ and $r\to\infty$.
\par
\vskip 2mm By  Lemma 4 and \eqref{eq3.15}, we get
\begin{eqnarray*}
\begin{aligned}
\quad &(j+1)T(r,F)\\
&\leq \overline{N}(r,F)+\overline{N}\left(r,\frac{1}{F-a+A}\right)+\overline{N}\left(r,\frac{1}{F-b+A}\right)+\sum_{i=1}^{j}\overline{N}\left(r,\frac{1}{F-d_{i}}\right)o(T(r,F))\\
&\leq T(r,F)+\sum_{i=1}^{j}\overline{N}\left(r,\frac{1}{F-d_{i}}\right)+o(T(r,F))\\
&\leq T(r,F)+\sum_{i=1}^{j}\overline{N}\left(r,\frac{1}{F-d_{i}}\right)+o(T(r,F)),
\end{aligned}
\end{eqnarray*}
as $r\to\infty$, possibly outside a set  $ E_{1}\cup E_{2}$  with $\int_{E_{2}}\rm dr<\infty$. This deduces
\begin{align}
\quad jT(r,F)\leq \sum_{i=1}^{j}\overline{N}\left(r,\frac{1}{F-d_{i}}\right)+o(T(r,F)),\label{eq3.16}
\end{align}
 as $r\not\in E_{1}\cup E_{2}\cup E_{3}$ and $r\to\infty$. According to $\{d_{1},d_{2},\ldots,d_{q}\}\subset\{d_{1},d_{2},\ldots,d_{j}\}$, \eqref{eq1}, \eqref{eq2}, \eqref{eq4}, \eqref{eq3.16},  and Lemma 3, we have
\begin{eqnarray*}
\begin{aligned}
&2T(r,F)=\bar{m}_{0,q}(r,F)+\sum_{d_{i}\in\hat{\mathbb{C}}}\hat{N}_{1}\left(r,\frac{1}{F-d_{i}}\right)+o(T(r,F))\\
&\leq \sum_{l=1}^{q}m\left(r,\frac{1}{F-d_{l}}\right)+m(r,F)+N(r,F)-\overline{N}(r,F)+\sum_{i=1}^{j}N\left(r,\frac{1}{F-d_{i}}\right)\\
&\quad-\sum_{i=1}^{j}\overline{N}\left(r,\frac{1}{F-d_{i}}\right)+o(T(r,F))=T(r,F)-\overline{N}(r,F)\\
&\quad+\sum_{l=1}^{q}m\left(r,\frac{1}{F-d_{i}}\right)+\sum_{i=1}^{j}N\left(r,\frac{1}{F-d_{i}}\right)-jT(r,F)\\
&\leq T(r,F)+\sum_{l=1}^{q}m\left(r,\frac{1}{F-d_{l}}\right)-\sum_{i=1}^{j}m\left(r,\frac{1}{F-d_{i}}\right)+o(T(r,F))\\
&\leq T(r,F)+o(T(r,F)),
\end{aligned}
\end{eqnarray*}
as $r\not\in E_{1}\cup E_{2}\cup E_{3}\subset(0,\infty)$ and $r\to\infty$. This  deduces that
\begin{align}
T(r,F)\leq o(T(r,F)),\label{eq3.17}
\end{align}
as $r\not\in E_{1}\cup E_{2}\cup E_{3}$ and $r\to\infty$. But it is impossible.
\par
\vskip 2mm This completes Theorem 4.

\section{The proof of Theorem 5 }
We prove Theorem 5 by contradiction. Suppose that $f\not\equiv \Delta_{c}^{n}f$. If $f$ is a non-constant polynomial, then $a$ and $b$ must be two finite values. By Lemma 2, we obtain $f\equiv \Delta_{c}^{n}f$, a contradiction. So $f$ is transcendental.  Since $f$ is a transcendental entire function, and that $f$ and $\Delta_{c}^{n}f$ share $a$ and $b$ IM, then  by the First Fundamental Theorem, Lemma 1 and Lemma 12  we get
\begin{eqnarray*}
\begin{aligned}
T(r,f)&\leq \overline{N}(r,\frac{1}{f-a})+\overline{N}(r,\frac{1}{f-b})+o(T(r,f))\\
&= \overline{N}(r,\frac{1}{\Delta_{c}^{n}f-a})+\overline{N}(r,\frac{1}{\Delta_{c}^{n}f-b})+o(T(r,f))\\
&\leq N(r,\frac{1}{f-\Delta_{c}^{n}f})+o(T(r,f))\\
&\leq T(r,f-\Delta_{c}^{n}f)+o(T(r,f))= m(r,f-\Delta_{c}^{n}f)+o(T(r,f))\\
&\leq m(r,f)+m(r,1-\frac{\Delta_{c}^{n}f}{f})+o(T(r,f))\\
&\leq T(r,f)+o(T(r,f)),
\end{aligned}
\end{eqnarray*}
as $r\to\infty$, possibly outside a set  $ E_{1}\subset(2,\infty)$  with $\int_{E_{1}}\rm dr<\infty$. It follows that
\begin{eqnarray}
T(r,f)=\overline{N}(r,\frac{1}{f-a})+\overline{N}(r,\frac{1}{f-b})+o(T(r,f)),\label{eq11.1}
\end{eqnarray}
as  $r\not\in E_{1}$ and $r\to\infty$.
\par
\vskip 2mm Suppose 
\begin{eqnarray}
F=\frac{f-a}{b-a}-d,\quad G=\frac{\Delta_{c}^{n}f-a}{b-a}-d,\label{eq11.2}
\end{eqnarray}
where $d$ is a non-constant small function of $f$. Because $f$ and $\Delta_{c}^{n}f$ share $a$ and $b$ IM, then from \eqref{eq11.2} we know that $F$ and $G$ share $-d$ and $1-d$ IM*. Let Let $v(r)\ll p\leq\infty$  be a large enough positive integer, let $i=1,2,3,\ldots,p$ be an integer, and let $h_{i}$ be $p$ distinct finite values. From Lemma 4, \eqref{eq11.1} and \eqref{eq11.2} we have
\begin{eqnarray*}
\begin{aligned}
&2T(r,F)\leq \overline{N}(r,F)+\overline{N}\left(r,\frac{1}{F}\right)+\overline{N}\left(r,\frac{1}{F-+d-1}\right)\notag\\
&+\overline{N}\left(r,\frac{1}{F-h_{i}}\right)+\varepsilon_{i}T(r,F)\leq T(r,F)+\overline{N}\left(r,\frac{1}{F-h_{i}}\right)+\varepsilon_{i}T(r,F)+o(T(r,f)),
\end{aligned}
\end{eqnarray*}
as $r\to\infty$, possibly outside a set  $ E_{1}\cup E_{2}\subset(2,\infty)$  with $\int_{E_{2}}\rm dr<\infty$. This deduces that 
\begin{align}
\lim_{r\to\infty}\frac{\overline{N}\left(r,\frac{1}{F-h_{i}}\right)}{T(r,F)}\geq1,\label{eq11.3}
\end{align}
 as $r\not\in E_{1}\cup E_{2}$ and $r\to\infty$, and where $\varepsilon_{i}>0$ is any positive number for every $i\in\{1,2,3,\ldots,p\}$.  Let $q=v(r)$ be a positive integer,  since $v(r)\sim\left(\log^{+}\frac{T(r,f)}{\log r}\right)^{20}=o(T(r,f))$, we have $\{1,2,3,\ldots,q\}\subset\{1,2,3,\ldots,p\}$ and $\{h_{1},h_{2},\ldots,h_{q}\}\subset\{h_{1},h_{2},\ldots,h_{p}\}$.  Suppose $\hat{\mathbb{C}}=\{h_{1},h_{2},\ldots,h_{j}\}\cup \{\infty\}$. According to $\{h_{1},h_{2},\ldots,h_{q}\}\subset\{h_{1},h_{2},\ldots,h_{p}\}$, \eqref{eq1}, \eqref{eq2}, \eqref{eq4} and Lemma 3
\begin{eqnarray*}
\begin{aligned}
&2T(r,F)=\bar{m}_{0,q}(r,F)+\sum_{h\in\hat{\mathbb{C}}}\hat{N}_{1}\left(r,\frac{1}{F-h_{i}}\right)+o(T(r,f))\\
&\leq m(r,F)+\sum_{i=1}^{q}m\left(r,\frac{1}{F-h_{i}}\right)+\sum_{i=1}^{p}N\left(r,\frac{1}{F-h_{i}}\right)-\sum_{i=1}^{p}\overline{N}\left(r,\frac{1}{F-h_{i}}\right)+o(T(r,f))\\
&\leq (p+1)T(r,F)+\sum_{i=1}^{q}m\left(r,\frac{1}{F-h_{i}}\right)-\sum_{i=1}^{p}m\left(r,\frac{1}{F-h_{i}}\right)-\sum_{i=1}^{p}\overline{N}\left(r,\frac{1}{F-h_{i}}\right)\\
&\quad+o(T(r,f))\leq (p+1)T(r,F)-\sum_{i=1}^{p}\overline{N}\left(r,\frac{1}{F-h_{i}}\right)++o(T(r,f)),
\end{aligned}
\end{eqnarray*}
as $r\not\in E_{1}\cup E_{2}\cup E_{3}$ and $r\to\infty,$ where and in what follows, $E_{4}\subset(0,\infty)$ is a set of logarithmic density $0$. This yields
\begin{align}
T(r,F)+\sum_{i=1}^{p}\overline{N}\left(r,\frac{1}{F-h_{i}}\right)\leq pT(r,F)+o(T(r,f)),\label{eq11.4}
\end{align}
 as $r\to\infty$  and  $r\not\in E_{1}\cup E_{2}\cup E_{3}$.   Divide $T(r,F)$ and we take the limit of the both sides of \eqref{eq11.4}, and combining \eqref{eq11.2} we can get
\begin{eqnarray*}
\begin{aligned}
1+p\leq 1+\sum_{i=1}^{p}\lim_{r\to\infty}\frac{\overline{N}\left(r,\frac{1}{F-h_{i}}\right)}{T(r,F)}\leq p,
\end{aligned}
\end{eqnarray*}
as $r\to\infty$  and  $r\not\in E_{1}\cup E_{2}\cup E_{3}$. This implies that $1\leq0$, a contradiction.
\par
\vskip 2mm This completes Theorem 5.

\section{The proof of Theorem 6}
\par
\vskip 2mm We prove Theorem 1 by contradiction. Suppose that $f_{c}\not\equiv f^{(k)}$. If $f$ is a non-constant polynomial, then $a$ and $b$ must be two finite values. Applying Lemma 2, we obtain $f_{c}\equiv f^{(k)}$, a contradiction. So $f$ is transcendental.  Since $f$ is a transcendental entire function, and that $f_{c}$ and $f^{(k)}$ share $a$ and $b$ IM, then  by the First Fundamental Theorem and Lemma 1 and Lemma 12,  we get
\begin{eqnarray*}
\begin{aligned}
T(r,f_{c})&\leq \overline{N}\left(r,\frac{1}{f_{c}-a}\right)+\overline{N}\left(r,\frac{1}{f_{c}-b}\right)+o(T(r,f))\\
&= \overline{N}\left(r,\frac{1}{f^{(k)}-a}\right)+\overline{N}\left(r,\frac{1}{f^{(k)}-b}\right)+o(T(r,f))\\
&\leq N\left(r,\frac{1}{f_{c}-f^{(k)}}\right)+o(T(r,f))\\
&\leq T(r,f_{c}-f^{(k)})+o(T(r,f))= m(r,f_{c}-f^{(k)})+o(T(r,f))\\
&\leq m(r,f_{c})+m\left(r,1-\frac{f^{(k)}}{f_{c}}\right)+o(T(r,f))\\
&\leq T(r,f_{c})+o(T(r,f)),
\end{aligned}
\end{eqnarray*}
as $r\to\infty$, possibly outside a set  $ E_{1}\subset(2,\infty)$  with $\int_{E_{1}}\rm dr<\infty$.
Thus we have
\begin{eqnarray}
T(r,f_{c})=\overline{N}\left(r,\frac{1}{f_{c}-a}\right)+\overline{N}\left(r,\frac{1}{f_{c}-b}\right)+o(T(r,f)), \label{eq12.1}
\end{eqnarray}
as  $r\not\in E_{1}$ and $r\to\infty$.
\par
\vskip 2mm Set
\begin{eqnarray}
\varphi=\frac{L(f_{c})(f_{c}-f^{(k)})}{(f_{c}-a)(f_{c}-b)},\\ \label{eq12.2}
\psi(z)=\frac{L(f^{(k)})(f_{c}-f^{(k)})}{(f^{(k)}-a)(f^{(k)}-b)}.\label{eq12.3}
\end{eqnarray}
\par
\vskip 2mm Note that  $f_{c}$ and $f^{(k)}$ share $a$ and $b$ IM.  We can know that $N(r,\varphi)=o(T(r,f))$ by \eqref{eq12.2}. By  Lemma 15,  we have
\begin{eqnarray*}
\begin{aligned}
&T(r,\varphi)=m(r,\varphi)=m\left(r,\frac{L(f_{c})(f_{c}-f^{(k)})}{(f_{c}-a)(f_{c}-b)}\right)+o(T(r,f))\notag\\
&=m\left(r,\frac{L(f_{c})f_{c}}{(f_{c}-a)(f_{c}-b)}\frac{f_{c}-f^{(k)}}{f_{c}}\right)+o(T(r,f))\\
&\leq m\left(r,\frac{L(f_{c})f_{c}}{(f_{c}-a)(f_{c}-b)}\right)+m\left(r,\frac{f_{c}-f^{(k)}}{f_{c}}\right)+o(T(r,f))\\
&\leq m\left(r,\frac{2L(f_{c})f_{c}}{(f_{c}-a)(f_{c}-b)}\right)+o(T(r,f))\\
&\leq m\left(r,\frac{2L(f_{c})f_{c}}{(f_{c}-a)(f_{c}-b)}-\frac{L(f_{c})(a+b)}{(f_{c}-a)(f_{c}-b)}\right)+m\left(r,\frac{L(f_{c})(a+b)}{(f_{c}-a)(f_{c}-b)}\right)+o(T(r,f))\\
&=m\left(r,\frac{L(f_{c})(2f_{c}-a-b)}{(f_{c}-a)(f_{c}-b)}\right)+m\left(r,\frac{L(f_{c})}{f_{c}-a}-\frac{L(f_{c})}{f_{c}-b}\right)+o(T(r,f))\\
&=m\left(r,\frac{L(f_{c})}{f_{c}-a}+\frac{L(f_{c})}{f_{c}-b}\right)+o(T(r,f))=o(T(r,f)),
\end{aligned}
\end{eqnarray*}
as  $r\not\in E_{1}$ and $r\to\infty$. This implies
\begin{align}
T(r,\varphi)=o(T(r,f)),\label{eq12.4}
\end{align}
as  $r\not\in E_{1}$ and $r\to\infty$.
\par
\vskip 2mm Let $d=a-j(a-b)$, where $ j\neq 0,1$ is a positive integer. Obviously, by  Lemma 15, we have
\begin{align}
 &m\left(r,\frac{1}{f_{c}}\right)=m\left(r,\frac{L(f_{c})}{\varphi(f_{c}-a)(f_{c}-b)}\frac{f_{c}-f^{(k)}}{f_{c}}\right)+o(T(r,f))\notag\\
&\leq m\left(r,\frac{L(f_{c})f_{c}}{(f_{c}-a)(f_{c}-b)}\right)+o(T(r,f))=o(T(r,f)),\label{eq12.5}
\end{align}
as  $r\not\in E_{1}$ and $r\to\infty$. And thus
\begin{align}
 m\left(r,\frac{1}{f_{c}-d}\right)&=m\left(r,\frac{L(f_{c})(f_{c}-f^{(k)})}{(\varphi (f_{c}-a)(f_{c}-b)(f_{c}-d)}\right)\notag\\
&\leq m\left(r,\frac{L(f_{c})f_{c}}{(f_{c}-a)(f_{c}-b)(f_{c}-d)}\right)\notag\\
&+ m\left(r,1-\frac{f^{(k)}}{f_{c}}\right)+o(T(r,f))=o(T(r,f)),\label{eq12.6}
\end{align}
as  $r\not\in E_{1}$ and $r\to\infty$.
\par
\vskip 2mm We discuss the following three cases.
\par
\vskip 2mm {\bf Case 1.}\quad $a$ and $b$ are finite constants. We can see from \eqref{eq3.1} and Theorem 3.4 that for  $a, b$  in $A$, we have
\begin{align}
 N\left(r,\frac{1}{f_{c}-a}\right)+N\left(r,\frac{1}{f_{c}-b}\right)&\leq T(r,f_{c})+N\left(r,\frac{1}{f_{c}^{(k)}}\right)\notag\\
 &+\varepsilon T(r,f_{c})+o(T(r,f)),\label{eq12.7}
\end{align}
as $r\not\in E_{1}\cup E_{2}\subset(0,\infty)$ and $r\to\infty,$ where and in what follows, $E_{2}\subset(0,\infty)$ is a set of logarithmic density $0$. It follows from the first main theorem of Nevanlinna and above that
\begin{eqnarray*}
\begin{aligned}
  &T(r,f_{c})\leq m\left(r,\frac{1}{f_{c}-a}\right)+m\left(r,\frac{1}{f_{c}-b}\right)+N\left(r,\frac{1}{f_{c}^{(k)}}\right)+\varepsilon T(r,f_{c})\notag\\
  &=m\left(r,\frac{1}{f_{c}-a}+\frac{1}{f_{c}-b}\right)+N\left(r,\frac{1}{f_{c}^{(k)}}\right)+\varepsilon T(r,f_{c})+o(T(r,f))\notag\\
  &\leq m\left(\frac{f_{c}^{(k)}}{f_{c}-a}+\frac{f_{c}^{(k)}}{f_{c}-b}\right)+m\left(r,\frac{1}{f_{c}^{(k)}}\right)\notag\\
  &+N\left(r,\frac{1}{f_{c}^{(k)}}\right)+\varepsilon T(r,f_{c})+o(T(r,f))\notag\\
  &\leq T(r,f_{c}^{(k)})+o(T(r,f))\leq T(r,f^{(k)})+\varepsilon T(r,f_{c})+o(T(r,f)),
\end{aligned}
\end{eqnarray*}
as $r\to\infty$  and  $r\not\in E_{1}\cup E_{2}$. It follows from above that
\begin{align}
T(r,f_{c})\leq T(r,f^{(k)})+\varepsilon T(r,f_{c})+o(T(r,f)),\label{eq12.8}
\end{align}
as $r\to\infty$  and  $r\not\in E_{1}\cup E_{2}$.
\par
\vskip 2mm Then by First Fundamental Theorem, Lemma 4, \eqref{eq12.1} and \eqref{eq12.8}, we have
\begin{eqnarray*}
\begin{aligned}
&2T(r,f_{c})\leq 2T(r,f^{(k)})+\varepsilon T(r,f_{c})+o(T(r,f))\\
&\leq\overline{N}\left(r,\frac{1}{f^{(k)}-a}\right)+\overline{N}\left(r,\frac{1}{f^{(k)}-b}\right)+\overline{N}\left(r,\frac{1}{f^{(k)}-d}\right)+\varepsilon T(r,f_{c})+o(T(r,f))\\
&\leq \overline{N}\left(r,\frac{1}{f-a}\right)+\overline{N}\left(r,\frac{1}{f-b}\right)+T\left(r,\frac{1}{f^{(k)}-d}\right)-m\left(r,\frac{1}{f^{(k)}-d}\right)\\
&+o(T(r,f))\leq T(r,f_{c})+T(r,f^{(k)})-m\left(r,\frac{1}{f^{(k)}-d}\right)+\varepsilon T(r,f_{c})+o(T(r,f))\\
&\leq 2T(r,f_{c})-m\left(r,\frac{1}{f^{(k)}-d}\right)+\varepsilon T(r,f_{c})+o(T(r,f)),
\end{aligned}
\end{eqnarray*}
as $r\to\infty$  and  $r\not\in E_{1}\cup E_{2}$. Thus
\begin{eqnarray}
m\left(r,\frac{1}{f^{(k)}-d}\right)\leq\varepsilon T(r,f_{c})+o(T(r,f)),\label{eq12.9}
\end{eqnarray}
as $r\to\infty$  and  $r\not\in E_{1}\cup E_{2}$. By the First Fundamental Theorem,  Lemma 13,   \eqref{eq12.5}-\eqref{eq12.6}, \eqref{eq12.8}-\eqref{eq12.9} and  that $f$ is a transcendental  entire function, we obtain
\begin{eqnarray*}
\begin{aligned}
&m\left(r,\frac{f_{c}-d}{f^{(k)}-d}\right)\leq m\left(r,\frac{f_{c}}{f^{(k)}-d})+m(r,\frac{d}{f^{(k)}-d}\right)+log2\\
&\leq T\left(r,\frac{f_{c}}{f^{(k)}-d}\right)-N\left(r,\frac{f_{c}}{f^{(k)}-d}\right)+o(T(r,f))=m\left(r,\frac{f^{(k)}-d}{f_{c}}\right)+N\left(r,\frac{f^{(k)}-d}{f_{c}}\right)\\
&-N\left(r,\frac{f_{c}}{f^{(k)}-d}\right)+o(T(r,f))\leq N\left(r,\frac{1}{f_{c}}\right)-N\left(r,\frac{1}{f^{(k)}-d}\right)+o(T(r,f))\\
&=T\left(r,\frac{1}{f_{c}}\right)-T\left(r,\frac{1}{f^{(k)}-d}\right)+o(T(r,f))\\
&=T(r,f_{c})-T(r,f^{(k)})+o(T(r,f))\leq\varepsilon T(r,f_{c})+o(T(r,f)),
\end{aligned}
\end{eqnarray*}
as $r\to\infty$  and  $r\not\in E_{1}\cup E_{2}$. This implies that
\begin{eqnarray}
m\left(r,\frac{f_{c}-d}{f^{(k)}-d}\right)\leq\varepsilon T(r,f_{c})+o(T(r,f)),\label{eq12.10}
\end{eqnarray}
as $r\to\infty$  and  $r\not\in E_{1}\cup E_{2}$.
By \eqref{eq3.3}, we have
\begin{eqnarray}
\psi=\left[\frac{a-d}{a-b}\frac{L(f^{(k)})}{f^{(k)}-a}-\frac{b-d}{a-b}\frac{L(f^{(k)})}{f^{(k)}-b}\right]\left[\frac{f-d}{f^{(k)}-d}-1\right].\label{eq12.11}
\end{eqnarray}
Since $f$ is a transcendental entire function, and that $f_{c}$ and $f^{(k)}$ share $a$ and $b$ IM, we know that  $N(r,\psi)=o(T(r,f))$. Then by \eqref{eq12.10}-\eqref{eq12.11} and  Lemma 15, we  get
\begin{align}
&T(r,\psi)=m(r,\psi)+o(T(r,f))\leq m\left(r,\frac{a-d}{a-b}\frac{L(f^{(k)})}{f^{(k)}-a}\right)\notag\\
&+m\left(r,\frac{b-d}{a-b}\frac{L(f^{(k)})}{f^{(k)}-b}\right)+m\left(r,\frac{f_{c}-d}{f^{(k)}-d}-1\right)+\varepsilon T(r,f)+o(T(r,f))\notag\\
&\leq\varepsilon T(r,f_{c})+o(T(r,f)),\label{eq12.12}
\end{align}
as $r\to\infty$  and  $r\not\in E_{1}\cup E_{2}$. Now let $z_{1}$ be a zero of $f_{c}-a$ and $f^{(k)}-a$ with multiplicities $m$ and $l$, respectively. Using Taylor series expansions,  and by calculating we get $l\varphi(z_{1})-m\psi(z_{1})=0$. Let
\begin{align}
G_{l,m}=l\varphi-m\psi,\label{eq12.13}
\end{align}
where $m$ and $l$ are positive integers. Next, we consider two subcases.
\par
\vskip 2mm {\bf Subcase 1.1.1.} \quad $G_{l,m}\equiv0$ for some positive integers $m$ and $l$. That is $l\varphi\equiv m\psi$. Then  we have
$$l \left(\frac{L(f_{c})}{f_{c}-a}-\frac{L(f_{c})}{f_{c}-b}\right)\equiv m \left(\frac{L(f^{(k)})}{f^{(k)}-a}-\frac{L(f^{(k)})}{f^{(k)}-b}\right),$$
which implies that
\begin{align}
   \left(\frac{f_{c}-a}{f_{c}-b}\right)^{l}\equiv A\left(\frac{f^{(k)}-a}{f^{(k)}-b}\right)^{m},\label{eq12.14}
\end{align}
where $A$ is a nonzero constant. By \eqref{eq12.14}, we obtain
\begin{eqnarray}
lT(r,f_{c})=mT(r,f^{(k)})+o(T(r,f))\leq mT(r,f)+o(T(r,f)).\label{eq12.15}
\end{eqnarray}
as $r\to\infty$  and  $r\not\in E_{1}\cup E_{2}$. If $l>m$, then we can obtain $T(r,f)=o(T(r,f))$, but it is impossible. If $m>l$,  we set $m=pl+t$, where $p>0$ and $0\leq t<l$ are two integers. As $f_{c}$ and $f^{(k)}$ share $a$ and $b$ IM, we can
see that $p\leq2$ and $t=0$, otherwise, we will get
\begin{eqnarray*}
\begin{aligned}
(pl+t)T(r,f^{(k)})=lT(r,f_{c})+o(T(r,f))\leq 2lT(r,f^{(k)})+o(T(r,f)),
\end{aligned}
\end{eqnarray*}
that is $T(r,f^{(k)})=o(T(r,f))$, a contradiction. If $m=2l$, then take $\varepsilon=\frac{1}{10}$ in \eqref{eq12.8}, we get
\begin{eqnarray*}
\begin{aligned}
2T(r,f^{(k)})&=T(r,f_{c})+o(T(r,f))\leq T(r,f^{(k)})+\frac{1}{10} T(r,f)+o(T(r,f))\\
&\leq T(r,f^{(k)})+\frac{1}{5} T(r,f^{(k)})+o(T(r,f)),
\end{aligned}
\end{eqnarray*}
and thus we have $T(r,f^{(k)})=o(T(r,f))$, s $r\to\infty$  and  $r\not\in E_{1}\cup E_{2}$. It is a contradiction. So $p=1$, and $m=l+t$. If $t>0$, rewrite \eqref{eq12.14} as
$$\left(\frac{f_{c}-a}{f_{c}-b}\right)^{l}\equiv A\left(\frac{f^{(k)}-a}{f^{(k)}-b}\right)^{l+t}.$$
Let $z_{1}$ be a simple zero of $f_{c}-a$ (or $f_{q,c}-b$), then the zeros(poles) of multiplicities of left hand side is $l$, but the zeros(poles) of multiplicities of right  hand side is  $l+t$. So  $f_{c}-a$ (or $f_{c}-b$) has no simple zero. Let $z_{1}$ be a  zero of $f_{c}-a$ (or $f_{c}-b$) with multiplicities $2$, and a zero of $f^{(k)}-a$(or $f^{(k)}-b$) with multiplicities $d$ then the zeros(poles) of multiplicities of left hand side is $2l$, but the zeros(poles) of multiplicities of right  hand side is  $d(l+t)$. It is easy to see that $d\leq2$. If $d=1$, we have $l=t$, a contradiction. If $d=2$, then by $t>0$, we can obtain a contradiction.  So  $f_{c}-a$ (or $f_{c}-b$) has no  zero with multiplicities $2$. Therefore, by \eqref{eq12.1}, we have
\begin{align}
T(r,f_{c})&=\overline{N}\left(r,\frac{1}{f_{c}-a}\right)+\overline{N}\left(r,\frac{1}{f_{c}-b}\right)+o(T(r,f))\notag\\
&=\overline{N}_{(3}\left(r,\frac{1}{f_{c}-a}\right)+\overline{N}_{(3}\left(r,\frac{1}{f_{c}-b}\right)+o(T(r,f))\notag\\
&\leq\frac{1}{3}(N\left(r,\frac{1}{f_{c}-a}\right)+N\left(r,\frac{1}{f_{c}-b})\right)+o(T(r,f))\notag\\
&\leq \frac{2}{3}T(r,f_{c})+o(T(r,f)),\label{eq12.16}
\end{align}
as $r\to\infty$  and  $r\not\in E_{1}\cup E_{2}$. It yields $T(r,f)=o(T(r,f))$, as $r\to\infty$  and  $r\not\in E_{1}\cup E_{2}$. This is a contradiction.  Thus we get $m=l$, in other words, $f_{c}$ and $f^{(k)}$ share $a$ and $b$ CM. It is easy to see that $\frac{Cb-a}{C-1}\neq a$ and $\frac{Cb-a}{C-1}\neq b$. It follows that $N(r,\frac{1}{f_{c}-\frac{Cb-a}{C-1}})=0$. Then by the Lemma 4,
\begin{eqnarray*}
\begin{aligned}
2T(r,f_{c})&\leq \overline{N}(r,f_{c})+\overline{N}\left(r,\frac{1}{f_{c}-a}\right)+\overline{N}\left(r,\frac{1}{f_{c}-b}\right)+\overline{N}\left(r,\frac{1}{f_{c}-\frac{Cb-a}{C-1}}\right)+\varepsilon_{1} T(r,f_{c})\\
&\leq \overline{N}\left(r,\frac{1}{f_{c}-a}\right)+\overline{N}\left(r,\frac{1}{f_{c}-b}\right)+\varepsilon_{1} T(r,f_{c}),
\end{aligned}
\end{eqnarray*}
as $r\to\infty$  and  $r\not\in E_{1}\cup E_{2}$. That is
\begin{eqnarray}
2T(r,f_{c})\leq \overline{N}\left(r,\frac{1}{f_{c}-a}\right)+\overline{N}\left(r,\frac{1}{f_{c}-b}\right)+\varepsilon_{1} T(r,f_{c}),\label{eq12.17}
\end{eqnarray}
as $r\to\infty$  and  $r\not\in E_{1}\cup E_{2}$. And it  follows from  \eqref{eq3.1} and \eqref{eq12.17} that $\varepsilon_{1}=\frac{1}{2}$ that $T(r,f_{c})=o(T(r,f))$, as $r\to\infty$  and  $r\not\in E_{1}\cup E_{2}$. It is a contradiction.
\par
\vskip 2mm {\bf Subcase 1.1.2.} \quad $G_{l,m}\not\equiv0$ for any positive integers $m$ and $l$. From the above discussion, we know that a zero of $f_{c}-a$ and $f^{(k)}-a$ (or a zero of $f_{c}-b$ and $f^{(k)}-b$ ) with multiplicities $m$ and $l$, must be the zero of $l\varphi-m\psi$.  So we have
\begin{align}
&\overline{N}_{(m,l)}\left(r,\frac{1}{f_{c}-a}\right)+\overline{N}_{(m,l)}\left(r,\frac{1}{f_{c}-b}\right)\leq \overline{N}\left(r,\frac{1}{l\varphi-m\psi}\right)+o(T(r,f))\notag\\
&\leq T(r,l\varphi-m\psi)+o(T(r,f))=\varepsilon T(r,f_{c})+o(T(r,f)),\label{eq12.18}
\end{align}
as $r\to\infty$  and  $r\not\in E_{1}\cup E_{2}$. Thus by \eqref{eq12.1} \eqref{eq12.8} and  \eqref{eq12.18}, we  get
\begin{align}
&T(r,f_{c})\leq \overline{N}\left(r,\frac{1}{f_{c}-a)}\right)+\overline{N}\left(r,\frac{1}{f_{c}-b}\right)+o(T(r,f))&\notag\\
&\leq N_{(1,1)}\left(r,\frac{1}{f_{c}-a}\right)+\overline{N}_{(2,1)}\left(r,\frac{1}{f_{c}-a}\right)+\overline{N}_{(3,1)}\left(r,\frac{1}{f_{c}-a}\right)\notag\\
&+\overline{N}_{(4,1)}\left(r,\frac{1}{f_{c}-a}\right)+\overline{N}_{(5}\left(r,\frac{1}{f_{c}-a}\right)\notag\\
&+N_{(1,1)}\left(r,\frac{1}{f_{c}-b}\right)+\overline{N}_{(2,1)}\left(r,\frac{1}{f_{c}-b}\right)+\overline{N}_{(3,1)}\left(r,\frac{1}{f_{c}-b}\right)\notag\\
&+\overline{N}_{(4,1)}\left(r,\frac{1}{f_{c}-b}\right)+\overline{N}_{(5}\left(r,\frac{1}{f_{c}-b}\right)+o(T(r,f))\notag\\
&\leq\frac{1}{5}\left(N(r,\frac{1}{f_{c}-a}\right)+N\left(r,\frac{1}{f_{c}-b}\right)+N\left(r,\frac{1}{f^{(k)}-a}\right)+N\left(r,\frac{1}{f^{(k)}-b})\right)\notag\\
&+8\varepsilon T(r,f)+o(T(r,f))\leq \frac{4}{5}T(r,f_{c})+8\varepsilon T(r,f)+o(T(r,f)),\label{eq12.19}
\end{align}
as $r\to\infty$  and  $r\not\in E_{1}\cup E_{2}$. From \eqref{eq12.19} we have $\frac{1}{5}T(r,f_{c})=4\varepsilon T(r,f)+o(T(r,f))$, as $r\to\infty$  and  $r\not\in E_{1}\cup E_{2}$. Take $\varepsilon=1/50$, and we can deduce $T(r,f_{c})=o(T(r,f))$, as $r\to\infty$  and  $r\not\in E_{1}\cup E_{2}$. It is a contradiction.

\par
\vskip 2mm {\bf Case 2}\quad $a$ and $b$ are two non-constant small functions. Let $p\leq\infty$ be a large enough positive integer, let $i=1,2,3,\ldots,p$ be an integer, and let $h_{i}$ be $p$ distinct finite values. From Lemma 4 and \eqref{eq12.1} we have
\begin{eqnarray*}
\begin{aligned}
&2T(r,f_{c})\leq \overline{N}(r,f_{c})+\overline{N}\left(r,\frac{1}{f_{c}-a}\right)+\overline{N}\left(r,\frac{1}{f_{c}-b}\right)\notag\\
&+\overline{N}\left(r,\frac{1}{f_{c}-h_{i}}\right)+\varepsilon_{i}T(r,f)\leq T(r,f_{c})+\overline{N}\left(r,\frac{1}{f_{c}-h_{i}}\right)+\varepsilon_{i}T(r,f)+o(T(r,f)),
\end{aligned}
\end{eqnarray*}
as $r\to\infty$, possibly outside a set  $ E_{1}\cup E_{3}\subset(2,\infty)$  with $\int_{E_{3}}\rm dr<\infty$. This deduces that 
\begin{align}
\lim_{r\to\infty}\frac{\overline{N}\left(r,\frac{1}{f_{c}-h_{i}}\right)}{T(r,f_{c})}\geq1,\label{eq12.20}
\end{align}
 as $r\not\in E_{1}\cup E_{2}$ and $r\to\infty$, and where $\varepsilon_{i}>0$ is any positive number for every $i\in\{1,2,3,\ldots,p\}$.  Let $q=v(r)$ be a positive integer,  since $v(r)\sim\left(\log^{+}\frac{T(r,f)}{\log r}\right)^{20}=o(T(r,f))$, we have $\{1,2,3,\ldots,q\}\subset\{1,2,3,\ldots,p\}$ and $\{h_{1},h_{2},\ldots,h_{q}\}\subset\{h_{1},h_{2},\ldots,h_{p}\}$.  Suppose $\hat{\mathbb{C}}=\{h_{1},h_{2},\ldots,h_{j}\}\cup \{\infty\}$. According to $\{h_{1},h_{2},\ldots,h_{q}\}\subset\{h_{1},h_{2},\ldots,h_{p}\}$, \eqref{eq1}, \eqref{eq2}, \eqref{eq4} and Lemma 3
\begin{eqnarray*}
\begin{aligned}
&2T(r,f_{c})=\bar{m}_{0,q}(r,f_{c})+\sum_{h\in\hat{\mathbb{C}}}\hat{N}_{0}\left(r,\frac{1}{f_{c}-h_{i}}\right)+o(T(r,f))\\
&\leq m(r,f_{c})+\sum_{i=1}^{q}m\left(r,\frac{1}{f_{c}-h_{i}}\right)+\sum_{i=1}^{p}N\left(r,\frac{1}{f_{c}-h_{i}}\right)-\sum_{i=1}^{p}\overline{N}\left(r,\frac{1}{f_{c}-h_{i}}\right)+o(T(r,f))\\
&\leq (p+1)T(r,f_{c})+\sum_{i=1}^{q}m\left(r,\frac{1}{f_{c}-h_{i}}\right)-\sum_{i=1}^{p}m\left(r,\frac{1}{f_{c}-h_{i}}\right)-\sum_{i=1}^{p}\overline{N}\left(r,\frac{1}{f_{c}-h_{i}}\right)\\
&\quad+o(T(r,f))\leq (p+1)T(r,f)-\sum_{i=1}^{p}\overline{N}\left(r,\frac{1}{f_{c}-h_{i}}\right)+\varepsilon T(r,f),
\end{aligned}
\end{eqnarray*}
as $r\not\in E_{1}\cup E_{3}\cup E_{4}$ and $r\to\infty,$ where and in what follows, $E_{4}\subset(0,\infty)$ is a set of logarithmic density $0$. This yields
\begin{align}
T(r,f_{c})+\sum_{i=1}^{p}\overline{N}\left(r,\frac{1}{f_{c}-h_{i}}\right)\leq pT(r,f_{c})+o(T(r,f)),\label{eq12.21}
\end{align}
 as $r\to\infty$  and  $r\not\in E_{1}\cup E_{3}\cup E_{4}$.   Divide $T(r,f_{c})$ and we take the limit of the both sides of \eqref{eq12.21}, and combining \eqref{eq12.20} we can get
\begin{eqnarray*}
\begin{aligned}
1+p\leq 1+\sum_{i=1}^{p}\lim_{r\to\infty}\frac{\overline{N}\left(r,\frac{1}{f_{c}-h_{i}}\right)}{T(r,f_{c})}\leq p,
\end{aligned}
\end{eqnarray*}
as $r\to\infty$  and  $r\not\in E_{1}\cup E_{3}\cup E_{4}$. This implies that $1\leq0$, a contradiction.
\par
\vskip 2mm{\bf Case 3}\quad $a$ is a non-constant small function, and $b$ is a finite value. Suppose $\hat{\mathbb{C}}=\{h_{1},h_{2},\ldots,h_{j}\}\cup \{b,\infty\}$. According to $\{h_{1},h_{2},\ldots,h_{q}\}\subset\{h_{1},h_{2},\ldots,h_{p}\}$, \eqref{eq1}, \eqref{eq2}, \eqref{eq4} and Lemma 3
\begin{eqnarray*}
\begin{aligned}
&2T(r,f_{c})=\bar{m}_{0,q}(r,f_{c})+\sum_{h\in\hat{\mathbb{C}}}\hat{N}_{0}\left(r,\frac{1}{f_{c}-h_{i}}\right)+o(T(r,f))\leq m(r,f_{c})\\
&\quad+m\left(r,\frac{1}{f_{c}-b}\right)+N\left(r,\frac{1}{f_{c}-b}\right)+\sum_{i=1}^{q}m\left(r,\frac{1}{f_{c}-h_{i}}\right)+\sum_{i=1}^{p}N\left(r,\frac{1}{f_{c}-h_{i}}\right)\\
&\quad-\sum_{i=1}^{p}\overline{N}\left(r,\frac{1}{f_{c}-h_{i}}\right)-\overline{N}\left(r,\frac{1}{f_{c}-b}\right)+o(T(r,f))\leq (p+2)T(r,f_{c})\\
&\quad+\sum_{i=1}^{q}m\left(r,\frac{1}{f_{c}-h_{i}}\right)-\sum_{i=1}^{p}m\left(r,\frac{1}{f_{c}-h_{i}}\right)-\sum_{i=1}^{p}\overline{N}\left(r,\frac{1}{f_{c}-h_{i}}\right)-\overline{N}\left(r,\frac{1}{f_{c}-b}\right)\\
&\quad+o(T(r,f))\leq (p+1)T(r,f_{c})-\sum_{i=1}^{p}\overline{N}\left(r,\frac{1}{f_{c}-h_{i}}\right)-\overline{N}\left(r,\frac{1}{f_{c}-b}\right)+o(T(r,f)),
\end{aligned}
\end{eqnarray*}
as $r\not\in E_{1}\cup E_{3}\cup E_{5}$ and $r\to\infty,$ where and in what follows, $E_{5}\subset(0,\infty)$ is a set of logarithmic density $0$. This yields
\begin{align}
\overline{N}\left(r,\frac{1}{f_{c}-b}\right)+\sum_{i=1}^{p}\overline{N}\left(r,\frac{1}{f_{c}-h_{i}}\right)\leq pT(r,f_{c})+o(T(r,f)),\label{eq12.22}
\end{align}
 as $r\to\infty$  and  $r\not\in E_{1}\cup E_{3}\cup E_{5}$.   Divide $T(r,f_{c})$ and we take the limit of the both sides of \eqref{eq12.22}, and combining \eqref{eq12.20} we can get
\begin{eqnarray*}
\begin{aligned}
\lim_{r\to\infty}\frac{\overline{N}\left(r,\frac{1}{f_{c}-b}\right)}{T(r,f_{c})}+p&\leq \lim_{r\to\infty}\frac{\overline{N}\left(r,\frac{1}{f_{c}-b}\right)}{T(r,f_{c})}+\sum_{i=1}^{p}\lim_{r\to\infty}\frac{\overline{N}\left(r,\frac{1}{f_{c}-h_{i}}\right)}{T(r,f_{c})}\leq p,
\end{aligned}
\end{eqnarray*}
as $r\to\infty$  and  $r\not\in E_{1}\cup E_{3}\cup E_{4}$. This and \eqref{eq12.1} and the assumption that $f_{c}$ and $f^{(k)}$ share $a$ and $b$ IM yield that
\begin{align}
T(r,f_{c})&=\overline{N}\left(r,\frac{1}{f_{c}-a}\right)+\overline{N}\left(r,\frac{1}{f_{c}-b}\right)+o(T(r,f))\notag\\
&\leq \overline{N}\left(r,\frac{1}{f_{c}-a}\right)+\varepsilon T(r,f)\leq T(r,f^{(k)})+\varepsilon T(r,f),\label{eq12.23}
\end{align}
 as $r\to\infty$  and  $r\not\in E_{1}\cup E_{2}\cup E_{5}$.   It follows from \eqref{eq12.23}
\begin{align}
T(r,f_{c})\leq T(r,f^{(k)})+\varepsilon T(r,f),\label{eq12.24}
\end{align}
as $r\to\infty$  and  $r\not\in E_{1}\cup E_{2}\cup E_{5}$. In this case, we can use the same method in {\bf Case 1} to obtain a contradiction.

\section{The proof of Theorem 7 }
\par
\vskip 2mm Suppose that $f$ is not a quasi-M$\ddot{o}$bius transformation of $g$, since $f$ and $g$ share $(a_{1},b_{1}), (a_{2},b_{2})$ CM and $(a_{3},b_{3}), (a_{4},b_{4}), (a_{5},b_{5})$ IM,  if $(a_{i},b_{i})$ are distinct pairs of values for $i\in\{1,2,3,4,5\}$, then by Theorem 5.1, we can obtain $f$ is  a quasi-M$\ddot{o}$bius transformation of $g$, but this is a contradiction.
\par
\vskip 2mm So, in the following, we assume that not all of  $a_{i}$ or $b_{i}$ are constants for $i\in\{1,2,3,4,5\}$. Without loss of generality, we assume that for $i\in\{1,2,3,4,5\}$, not all of $a_{i}$  are constants. Let $p\leq\infty$ be a large enough positive integer, let $k=1,2,3,\ldots,p$ be an integer, and let $c_{k}\not\equiv a_{1}, a_{2}, a_{3}, a_{4}, a_{5}$ be $p$ distinct finite values.  By Lemma 4 and Lemma 18, we get
\begin{eqnarray*}
\begin{aligned}
4T(r,f)&\leq \sum_{i=1}^{5}\overline{N}\left(r,\frac{1}{f-a_{i}}\right)+\overline{N}\left(r,\frac{1}{f-c_{k}}\right)+S^{*}(r,f)\\
&=3T(r,f)+\overline{N}\left(r,\frac{1}{f-c_{k}}\right)+S^{*}(r,f),
\end{aligned}
\end{eqnarray*}
as $r\to\infty$, possibly outside a set  $ E_{1}\subset(2,\infty)$  with $\int_{E_{1}}\rm dr<\infty$. It follows from above that
\begin{align}
T(r,f)= \overline{N}\left(r,\frac{1}{f-c_{k}}\right)+S^{*}(r,f),\label{eq13.1}
\end{align}
as  $r\not\in E_{1}$ and $r\to\infty$. 
By \eqref{eq3.1} , we get
\begin{align}
\lim_{r\to\infty}\frac{\overline{N}\left(r,\frac{1}{f-c_{k}}\right)}{T(r,f)}\geq1,\label{eq13.2}
\end{align}
as $r\not\in E_{1}\cup E_{2}$ and $r\to\infty$.
Let $q=v(r)$ be a positive integer,  and as $v(r)\sim(\log^{+}\frac{T(r,f)}{\log r})^{20}=o(T(r,f))$, we have $\{1,2,3,\ldots,q\}\subset\{1,2,3,\ldots,p\}$ and $\{c_{1},c_{2},\ldots,c_{q}\}\subset\{c_{1},c_{2},\ldots,c_{p}\}$.  Suppose $\hat{\mathbb{C}}=\{c_{1},c_{2},\ldots,c_{p}\}\cup \{\infty\}$. According to $\{c_{1},c_{2},\ldots,c_{q}\}\subset\{c_{1},c_{2},\ldots,c_{p}\}$. We discuss the following five cases.
\par
\vskip 2mm {\bf Case 1}\quad $a_{1}$ is a constant and $ a_{2}, a_{3}, a_{4}, a_{5}$ are non-constant small functions. Suppose $\hat{\mathbb{C}}=\{c_{1},c_{2},\ldots,c_{j}\}\cup \{a_{1}\}$. According to $\{c_{1},c_{2},\ldots,c_{q}\}\subset\{c_{1},c_{2},\ldots,c_{p}\}$, \eqref{eq1}, \eqref{eq2}, \eqref{eq4} and Lemma 3, we have
\begin{eqnarray*}
\begin{aligned}
&2T(r,f)=\bar{m}_{0,q}(r,f)+\sum_{c\in\widehat{\mathbb{C}}}\hat{N}_{1}(r,\frac{1}{f-c})+S(r,f)\leq m(r,a_{1},f)+N(r,a_{1},f)\\
&\leq +\sum_{k=1}^{q}m(r,c_{k},f)+\sum_{k=1}^{p}N(r,c_{k},f)-\sum_{k=1}^{p}\overline{N}(r,c_{k},f)-\overline{N}(r,a_{1},f)+S(r,f)\\
&\leq (p+1)T(r,f)+\sum_{k=1}^{q}m(r,c_{k},f)-\sum_{k=1}^{p}m(r,c_{k},f)-\sum_{k=1}^{p}\overline{N}(r,c_{k},f)-\overline{N}(r,a_{1},f)\\
&+S(r,f)\leq (p+1)T(r,f)-\sum_{k=1}^{p}\overline{N}(r,c_{k},f)-\overline{N}(r,a_{1},f)+S(r,f),
\end{aligned}
\end{eqnarray*}
as $r\not\in E_{1}\cup E_{2}\subset(0,\infty)$ and $r\to\infty,$ where and in what follows, $E_{2}\subset(0,\infty)$ is a set of logarithmic density $0$. This yields 
\begin{align}
T(r,f)+\sum_{k=1}^{p}\overline{N}(r,c_{k},f)+\overline{N}(r,a_{1},f)=pT(r,f),\label{eq13.3}
\end{align}
as $r\to\infty$  and  $r\not\in E_{1}\cup E_{2}$. Divide $T(r,f)$ and we take the infimum of the both sides of \eqref{eq13.3}, and combining \eqref{eq13.2}, we can get
 \begin{align}
&1+p\leq 1+\lim_{r\to\infty}\frac{\sum_{k=1}^{p}\overline{N}\left(r,\frac{1}{f-c_{k}}\right)}{T(r,f)}\leq p,\label{eq13.4}
\end{align}
  as $r\to\infty$  and  $r\not\in E_{1}\cup E_{2}$. From \eqref{eq3.4} we have $1\leq0$, but it is impossible.
\par
\vskip 2mm{\bf Case 2}\quad $a_{1}$ and $a_{2}$ are two constants and $a_{3}, a_{4}, a_{5}$ are non-constant small functions. Suppose $\hat{\mathbb{C}}=\{c_{1},c_{2},\ldots,c_{j}\}\cup \{a_{1}, a_{2}\}$. From $\{c_{1},c_{2},\ldots,c_{q}\}\subset\{c_{1},c_{2},\ldots,c_{p}\}$, \eqref{eq1}, \eqref{eq2}, \eqref{eq4} and Lemma 3, we have
\begin{eqnarray*}
\begin{aligned}
&2T(r,f)=\bar{m}_{0,q}(r,f)+\sum_{c\in\hat{\mathbb{C}}}\hat{N}_{1}(r,\frac{1}{f-c})+S(r,f)\leq m(r,a_{1},f)+m(r,a_{2},f)\\
&+N(r,a_{1},f)+N(r,a_{1},f)+\sum_{k=1}^{q}m(r,c_{k},f)+\sum_{k=1}^{p}N(r,c_{k},f)-\overline{N}(r,a_{1},f)\\
&-\overline{N}(r,a_{2},f)-\sum_{k=1}^{p}\overline{N}(r,c_{k},f)+S(r,f)\leq (p+2)T(r,f)+\sum_{k=1}^{q}m(r,c_{k},f)\\
&-\sum_{k=1}^{p}m(r,c_{k},f)-\sum_{k=1}^{p}\overline{N}(r,c_{k},f)-\overline{N}(r,a_{1},f)-\overline{N}(r,a_{2},f)+S(r,f)\\
&\leq (p+2)T(r,f)-\sum_{k=1}^{p}\overline{N}(r,c_{k},f)-\overline{N}(r,a_{1},f)-\overline{N}(r,a_{2},f)+S(r,f),
\end{aligned}
\end{eqnarray*}
as $r\not\in E_{1}\cup E_{3}\subset(0,\infty)$ and $r\to\infty,$ where and in what follows, $E_{3}\subset(0,\infty)$ is a set of logarithmic density $0$. This obtains
\begin{align}
\overline{N}(r,a_{1},f)+\overline{N}(r,a_{2},f)+\sum_{k=1}^{p}\overline{N}(r,c_{k},f)\leq p(r,f)+S(r,f),\label{eq13.5}
\end{align}
as $r\to\infty$  and  $r\not\in E_{1}\cup E_{3}$.  Divide $T(r,f)$ and we take the limit of the both sides of \eqref{eq13.3}, and combining \eqref{eq13.5}, we can get
 \begin{align}
&\lim_{r\to\infty}\frac{\overline{N}(r,a_{1},f)+\overline{N}(r,a_{2},f)}{T(r,f)}+p\leq \lim_{r\to\infty}\frac{\overline{N}(r,a_{1},f)+\overline{N}(r,a_{2},f)}{T(r,f)}\notag\\
&+\lim_{r\to\infty}\frac{\sum_{k=1}^{p}\overline{N}\left(r,\frac{1}{f-c_{k}}\right)}{T(r,f)}\leq p,\label{eq13.6}
\end{align}
  as $r\to\infty$  and  $r\not\in E_{1}\cup E_{3}$. \eqref{eq13.6} yields
\begin{eqnarray*}
\begin{aligned}
\lim_{r\to\infty}\frac{\overline{N}(r,a_{1},f)+\overline{N}(r,a_{2},f)}{T(r,f)}\leq0,
\end{aligned}
\end{eqnarray*}
that is
\begin{align}
\overline{N}(r,a_{1},f)+\overline{N}(r,a_{2},f)\leq S^{*}(r,f),\label{eq13.7}
\end{align}
Since $\varepsilon>0$ is arbitrary in $S^{*}(r,f)$, we can obtain from Lemma 18 and \eqref{eq13.7} that $T(r,f)=S(r,f)$, as $r\to\infty$  and  $r\not\in E_{1}\cup E_{3}$. This is impossible.

{\bf Case 3}\quad $a_{1}, a_{2}$ and $a_{3}$ are three constants and $a_{4}, a_{5}$ are non-constant small functions.  Suppose $\hat{\mathbb{C}}=\{c_{1},c_{2},\ldots,c_{j}\}\cup \{a_{1}, a_{2}, a_{3}\}$. From $\{c_{1},c_{2},\ldots,c_{q}\}\subset\{c_{1},c_{2},\ldots,c_{p}\}$, \eqref{eq1}, \eqref{eq2}, \eqref{eq4} and Lemma 3, we have
\begin{eqnarray*}
\begin{aligned}
&2T(r,f)=\bar{m}_{0,q}(r,f)+\sum_{c\in\hat{\mathbb{C}}}\hat{N}_{1}(r,\frac{1}{f-c})+S(r,f)\leq m(r,a_{1},f)+m(r,a_{2},f)\\
&+m(r,a_{3},f)+N(r,a_{1},f)+N(r,a_{2},f)+N(r,a_{3},f)+\sum_{k=1}^{q}m(r,c_{k},f)\\
&+\sum_{k=1}^{p}N(r,c_{k},f)-\sum_{k=1}^{p}\overline{N}(r,c_{k},f)+S(r,f)\leq (p+3)(r,f)+\sum_{k=1}^{q}m(r,c_{k},f)\\
&-\sum_{k=1}^{p}m(r,c_{k},f)-\overline{N}(r,a_{1},f)-\overline{N}(r,a_{2},f)-\overline{N}(r,a_{3},f)+(r,f)\\
&\leq (p+3)T(r,f)-\sum_{k=1}^{p}\overline{N}(r,c_{k},f)-\overline{N}(r,a_{1},f)-\overline{N}(r,a_{2},f)-\overline{N}(r,a_{3},f)+S(r,f)
\end{aligned}
\end{eqnarray*}
as $r\not\in E_{1}\cup E_{4}\subset(0,\infty)$ and $r\to\infty,$ where and in what follows, $E_{4}\subset(0,\infty)$ is a set of logarithmic density $0$. This obtains
\begin{align}
\overline{N}(r,a_{1},f)+\overline{N}(r,a_{2},f)+\overline{N}(r,a_{3},f)+\sum_{k=1}^{p}\overline{N}(r,c_{k},f)\leq (p+1)T(r,f)+S^{*}(r,f),\label{eq13.8}
\end{align}
as $r\to\infty$  and  $r\not\in E_{1}\cup E_{4}$. Divide $T(r,f)$ and we take the limit of the both sides of \eqref{eq13.3}, and combining \eqref{eq13.8}, we can get
 \begin{align}
&\lim_{r\to\infty}\frac{\overline{N}(r,a_{1},f)+\overline{N}(r,a_{2},f)+\overline{N}(r,a_{3},f)}{T(r,f)}+p\leq \lim_{r\to\infty}\frac{\sum_{k=1}^{p}\overline{N}\left(r,\frac{1}{f-c_{k}}\right)}{T(r,f)}\notag\\
&+\lim_{r\to\infty}\frac{\overline{N}(r,a_{1},f)+\overline{N}(r,a_{2},f)+\overline{N}(r,a_{3},f)}{T(r,f)}\leq (p+1),\label{eq13.9}
\end{align}
  as $r\to\infty$  and  $r\not\in E_{1}\cup E_{4}$. \eqref{eq13.9} yields
\begin{eqnarray*}
\begin{aligned}
\lim_{r\to\infty}\frac{\overline{N}(r,a_{1},f)+\overline{N}(r,a_{2},f)+\overline{N}(r,a_{3},f)}{T(r,f)}\leq1,
\end{aligned}
\end{eqnarray*}
that is
\begin{align}
\overline{N}(r,a_{1},f)+\overline{N}(r,a_{2},f)+\overline{N}(r,a_{3},f)\leq T(r,f)+S^{*}(r,f),\label{eq13.10}.
\end{align}
as $r\to\infty$  and  $r\not\in E_{1}\cup E_{4}$. Since $\varepsilon>0$ is arbitrary in $S^{*}(r,f)$, we can obtain from the third equality of Lemma 17 and \eqref{eq13.10} that $T(r,f)=S(r,f)$, as $r\to\infty$  and  $r\not\in E_{1}\cup E_{4}$. This is a contradiction.
\par
\vskip 2mm{\bf Case 4}\quad $a_{1}, a_{2},a_{3}$ and $a_{4}$ are four constants and $a_{5}$ is a non-constant small function.  We set
\begin{align}
F=f-a_{5},\quad G=g-b_{5}\label{eq13.11}.
\end{align}
If one of $a_{1}, a_{2},a_{3}, a_{4}$ is $\infty$, says $a_{1}=\infty$,  we set
then we know from \eqref{eq13.11} that $F$ and $G$ share $(\infty, b_{1}-b_{5})$ and  $(a_{2}-a_{5}, b_{2}-b_{5})$ CM*, and share $(a_{3}-a_{5}, b_{3}-b_{5})$,   $(a_{4}-a_{5}, b_{4}-b_{5})$ and $(0, 0)$ IM*. Another thing is that $a_{i}-a_{5}$ are non-constant small functions for $i\in\{2,3,4\}$ of $F$. Similar to the proof of {\bf Case 2}, we can obtain $T(r,f)=S(r,f)$, a contradiction. If none of $a_{1}, a_{2},a_{3}, a_{4}$ is $\infty$, then $F$ and $G$ share $(a_{1}-a_{5}, b_{1}-b_{5})$ and  $(a_{2}-a_{5}, b_{2}-b_{5})$ CM*, and share $(a_{3}-a_{5}, b_{3}-b_{5})$,   $(a_{4}-a_{5}, b_{4}-b_{5})$ and $(0, 0)$ IM*. Another thing is that $a_{i}-a_{5}$ are non-constant small functions for $i\in\{1,2,3,4\}$ of $F$. Similar to the proof of {\bf Case 1}, we can obtain $T(r,f)=S(r,f)$, a contradiction.
\par
\vskip 2mm{\bf Case 5}\quad $a_{1}, a_{2}, a_{3}, a_{4}, a_{5}$ are non-constant small functions. Suppose we have $\hat{\mathbb{C}}=\{c_{1},c_{2},\ldots,c_{j}\}$. From $\{c_{1},c_{2},\ldots,c_{q}\}\subset\{c_{1},c_{2},\ldots,c_{p}\}$, \eqref{eq1}, \eqref{eq2}, \eqref{eq4} and Lemma 3, we have
\begin{eqnarray*}
\begin{aligned}
&2T(r,f)=\bar{m}_{0,q}(r,f)+\sum_{c\in\hat{\mathbb{C}}}\hat{N}_{1}(r,\frac{1}{f-c})+S(r,f)\leq \sum_{k=1}^{q}m(r,c_{k},f)+\sum_{k=1}^{p}N(r,c_{k},f)\\
&-\sum_{k=1}^{p}\overline{N}(r,c_{k},f)+S(r,f)\leq pT(r,f)+\sum_{k=1}^{q}m(r,c_{k},f)-\sum_{k=1}^{p}m(r,c_{k},f)\\
&-\sum_{k=1}^{p}\overline{N}(r,c_{k},f)+S(r,f)\leq pT(r,f)-\sum_{k=1}^{p}\overline{N}(r,c_{k},f)+ S(r,f),
\end{aligned}
\end{eqnarray*}
as $r\not\in E_{1}\cup E_{5}\subset(0,\infty)$ and $r\to\infty,$ where and in what follows, $E_{5}\subset(0,\infty)$ is a set of logarithmic density $0$. This obtains
\begin{align}
2T(r,f)+\sum_{k=1}^{p}\overline{N}(r,c_{k},f)\leq pT(r,f)+ S(r,f),\label{eq13.12}
\end{align}
as $r\to\infty$  and  $r\not\in E_{1}\cup E_{5}$.  Divide $T(r,f)$ and we take the limit of the both sides of \eqref{eq13.3}, and combining \eqref{eq13.12}, we can get
 \begin{align}
&2+p\leq 2+\lim_{r\to\infty}\frac{\sum_{k=1}^{p}\overline{N}\left(r,\frac{1}{f-c_{k}}\right)}{T(r,f)}\leq p,\label{eq13.13}
\end{align}
  as $r\to\infty$  and  $r\not\in E_{1}\cup E_{2}$. From \eqref{eq3.13} we have $2\leq0$, but it is impossible.
\par
\vskip 2mm This completes Theorem 7.

\section{The proof of Theorem 8 }
Suppose $f^{(k)}\not\equiv g^{(k)}$. If $f$ is rational, then $T(r,f)=O(\log r)$. In this case,  three distinct small functions $a_{1}$, $a_{2}$, and $a_{3}$ of $f$ must be three distinct finite values. Thus, by Lemma 9, we can obtain $f^{(k)}\equiv g^{(k)}$, a contradiction. And hence $f$ is transcendental. By Lemma 4 and (ii) of Lemma 5, we can obtain that for  $i=1,2,3,\ldots,p$,
 $d_{i}\not\equiv a_{1}, a_{2}, a_{3}, \infty$ are distinct values,
\begin{eqnarray*}
\begin{aligned}
(p+2)T(r,f^{(k)})&\leq \overline{N}(r,f^{(k)})+\overline{N}(r,\frac{1}{f^{(k)}-a_{1}})+\overline{N}(r,\frac{1}{f^{(k)}-a_{2}})+\overline{N}(r,\frac{1}{f^{(k)}-a_{3}})\notag\\
&+\sum_{i=1}^{p}\overline{N}(r,\frac{1}{f^{(k)}-d_{i}})+\varepsilon T(r,f^{(k)})+o(T(r,f^{(k)}))\notag\\
&\leq (2+\varepsilon_{1}+\varepsilon)T(r,f^{(k)})+\sum_{i=1}^{p}\overline{N}(r,\frac{1}{f^{(k)}-d_{i}})+o(T(r,f^{(k)}))
\end{aligned}
\end{eqnarray*}
that is
\begin{align}
pT(r,f^{(k)})\leq \sum_{i=1}^{p}\overline{N}(r,\frac{1}{f^{(k)}-d_{i}})+\varepsilon_{3} T(r,f^{(k)})+o(T(r,f^{(k)})),
\end{align}
for all $\varepsilon_{1}+\varepsilon=\varepsilon_{3}>0$ outside a set $E_{1}$ with $\int_{E_{1}}d r<\infty$. Let $q=v(r)$ be a positive integer, and  as $v(r)\sim(\log^{+}\frac{T(r,f^{(k)})}{\log r})^{20}=o(T(r,f^{(k)}))$, we discuss the following four cases.
\par
\vskip 2mm Next, we discuss four cases.
\par
\vskip 2mm {\bf Case 1} $a_{1}=0, a_{2}=1, a_{3}=b$ are three distinct finite values. Suppose we have
\begin{align}
F=f^{(k)}-A,\quad G=g^{(k)}-A,\label{eq14.1}
\end{align}
where $A$ is a non-constant rational function. Since $f^{(k)}$ and $g^{(k)}$ share $0,1,b$ IM* and $\infty$ IM*, we know that $F$ and $G$ share $-A,1-A,b-A$ IM* and $\infty$ IM*. Let $v(r)\ll p\leq\infty$ be a large enough positive integer, let $j=1,2,3,\ldots,p$ be an integer,  and let  $d_{j}\not\equiv 0,1,b,\infty$ are distinct values.   Applying Lemma 4 to $F$ for every $j\in\{1,2,3,\ldots,p\}$, we can obtain 
\begin{eqnarray*}
\begin{aligned}
&(p+2)T(r,F)\leq \overline{N}\left(r,F\right)+\overline{N}\left(r,\frac{1}{F+A}\right)+\overline{N}\left(r,\frac{1}{F-1+A}\right)+\overline{N}\left(r,\frac{1}{F-b+A}\right)\\
&\quad+\sum_{j=1}^{p}\overline{N}\left(r,\frac{1}{F-d_{j}}\right)+S(r,F)\leq 4T(r,F)+\sum_{j=1}^{p}\overline{N}\left(r,\frac{1}{F-d_{j}}\right)\\
&-m(r,F)-m\left(r,\frac{1}{F+A}\right)-m\left(r,\frac{1}{F-1+A}\right)-m\left(r,\frac{1}{F-b+A}\right)+S(r,F),
\end{aligned}
\end{eqnarray*}
as $r\to\infty$, possibly outside a set  $ E_{1}\subset(2,\infty)$  with $\int_{E_{1}}\rm dr<\infty$. This deduces
\begin{align}
&(p-2)T(r,F)+m(r,F)+m\left(r,\frac{1}{F+A}\right)+m\left(r,\frac{1}{F-1+A}\right)\notag\\
&+m\left(r,\frac{1}{F-b+A}\right)\leq \sum_{j=1}^{p}\overline{N}\left(r,\frac{1}{F-d_{j}}\right)+S(r,F),\label{eq14.2}
\end{align}
as $r\to\infty$ and $r\not\in E_{1}$. 
\par
\vskip 2mm Let $q=v(r)$ be a positive integer,  and as $v(r)\sim(\log^{+}\frac{T(r,f)}{\log r})^{20}=o(T(r,f))$, we have $\{1,2,3,\ldots,q\}\subset\{1,2,3,\ldots,p\}$ and $\{d_{1},d_{2},\ldots,d_{q}\}\subset\{d_{1},d_{2},\ldots,d_{p}\}$. Suppose $\hat{\mathbb{C}}=\{d_{1},d_{2},\ldots,d_{p}\}\cup \{a_{3}, a_{4}\}$. According to $\{d_{1},d_{2},\ldots,d_{q}\}\subset\{d_{1},d_{2},\ldots,d_{p}\}$, \eqref{eq1}, \eqref{eq2}, \eqref{eq4}, \eqref{eq14.2}, (v) of Lemma 5 and Lemma 3, we have
\begin{eqnarray*}
\begin{aligned}
&2T(r,F)=\bar{m}_{0,q}(r,F)+\sum_{d\in\hat{\mathbb{C}}}\hat{N}_{1}\left(r,\frac{1}{F-d_{j}}\right)+S(r,F)\leq m\left(r,F\right)+N\left(r,F\right)\\
&+\sum_{j=1}^{q}m\left(r,\frac{1}{F-d_{j}}\right)
-\overline{N}\left(r,F\right)+\sum_{j=1}^{p}N\left(r,\frac{1}{F-d_{j}}\right)
-\sum_{j=1}^{p}\overline{N}\left(r,\frac{1}{F-d_{j}}\right)+S(r,F)\\
&=(p+1)T(r,F)+\sum_{j=1}^{q}m\left(r,\frac{1}{F-d_{j}}\right)-\sum_{j=1}^{p}m\left(r,\frac{1}{F-d_{j}}\right)-(p-2)T(r,F)\\
&-\overline{N}\left(r,f\right)-m(r,F)-m\left(r,\frac{1}{F+A}\right)-m\left(r,\frac{1}{F-1+A}\right)-m\left(r,\frac{1}{F-b+A}\right)+S(r,F)\\
&\leq 2T(r,F)-m\left(r,\frac{1}{F+A}\right)-m\left(r,\frac{1}{F-1+A}\right)-m\left(r,\frac{1}{F-b+A}\right)+S(r,F)
\end{aligned}
\end{eqnarray*}
as $r\not\in E_{1}\cup E_{2}$ and $r\to\infty,$ where and in what follows, $E_{2}\subset(0,\infty)$ is a set of logarithmic density $0$. This and \eqref{eq12} deduce that
\begin{align}
m\left(r,\frac{1}{f^{(k)}}\right)+m\left(r,\frac{1}{f^{(k)}-1}\right)+m\left(r,\frac{1}{f^{(k)}-b}\right)\leq S(r,f),\label{eq14.3}
\end{align}
as  $r\not\in E_{1}\cup E_{2}$ and $r\to\infty$. Similarly, we also have
\begin{align}
m\left(r,\frac{1}{g^{(k)}}\right)+m\left(r,\frac{1}{g^{(k)}-1}\right)+m\left(r,\frac{1}{g^{(k)}-b}\right)\leq S(r,f)\leq S(r,g)=S(r,f),\label{eq14.4}
\end{align}
as $r\not\in E_{3}\cup E_{4}$ and $r\to\infty,$ where and in what follows, $ E_{3}\subset(2,\infty)$  with $\int_{E_{3}}\rm dr<\infty$ and $E_{4}\subset(0,\infty)$ is a set of logarithmic density $0$.
\par
\vskip 2mm Let
\begin{align}
F_{1}=\frac{1}{f^{(k)}},\quad G_{1}=\frac{1}{g^{(k)}}.\label{eq14.5}
\end{align}
Because  $f^{(k)}$ and $g^{(k)}$ share $0,1,b$ IM* and $\infty$ IM*, we know that $F_{1}$ and $G_{1}$ share $c=\frac{1}{b},1, \infty$ IM*, and $0$ CM*. Suppose we have the following two functions,
\begin{align}
\varphi=\frac{F'_{1}(F_{1}-G_{1})}{F_{1}(F_{1}-1)(F_{1}-c)},\label{eq14.6}
\end{align}
\begin{align}
\psi=\frac{G'_{1}(F_{1}-G_{1})}{G_{1}(G_{1}-1)(G_{1}-c)}.\label{eq14.7}
\end{align}
It is easy to see that $\varphi\not\equiv0$ and $\psi\not\equiv0$ since $F_{1}\not\equiv G_{1}$. On the one hand,
\begin{align}
m(r,\varphi)&=m\left(r,\frac{F'_{1}(F_{1}-G_{1})}{F_{1}(F_{1}-1)(F_{1}-c)}\right)\notag\\
&\leq m\left(r,\frac{F'_{1}F_{1}}{F_{1}(F_{1}-1)(F_{1}-c)}\right)+m\left(r,\frac{F'_{1}}{F_{1}(F_{1}-1)(F_{1}-c)}\right)+m(r,G_{1})\notag\\
&\leq m\left(r,\frac{1}{g}\right)+S(r,f)=S(r,f),\label{eq14.8}
\end{align}
as $r\not\in E_{1}\cup E_{2}\cup E_{3}\cup E_{4}$ and $r\to\infty$. Similarly, we have
\begin{align}
m(r,\psi)\leq m\left(r,\frac{1}{f}\right)+S(r,f)=S(r,f),\label{eq14.9}
\end{align}
as $r\not\in E_{1}\cup E_{2}\cup E_{3}\cup E_{4}$ and $r\to\infty$. Next, we analyze the poles of $\varphi$ and $\psi$. As  $F_{1}$ and $G_{1}$ share $1,c,\infty$ IM* and $0$ CM*, so we can see from \eqref{eq14.6} and \eqref{eq14.7} that
\begin{align}
N(r,\varphi)\leq \sum_{n=2}N_{(1,n)}(r,G_{1})-\sum_{n=2}\overline{N}_{(1,n)}(r,G_{1}),\label{eq14.10}
\end{align}
\begin{align}
N(r,\psi)\leq \sum_{m=2}N_{(m,1)}(r,F_{1})-\sum_{m=2}\overline{N}_{(m,1)}(r,F_{1}),\label{eq14.11}
\end{align}
where $m$ and $n$ are two positive integers.
\par
\vskip 2mm Note that the poles of $F$ with multiplicities $m\geq2$ and the poles of $G$ with multiplicities $1$ are the zero points of $\varphi$, so we get that
\begin{align}
&\sum_{m=2}N_{(m,1)}(r,F_{1})\leq N\left(r,\frac{1}{\varphi}\right)\leq T(r,\varphi)=N(r,\varphi)+S(r,f)\notag\\
&\leq \sum_{n=2}N_{(1,n)}(r,G_{1})-\sum_{n=2}\overline{N}_{(1,n)}(r,G_{1})+S(r,f),\label{eq14.12}
\end{align}
as $r\not\in E_{1}\cup E_{2}\cup E_{3}\cup E_{4}$ and $r\to\infty$. Similarly,
\begin{align}
&\sum_{n=2}N_{(1,n)}(r,G_{1})\leq N\left(r,\frac{1}{\psi}\right)\leq T(r,\psi)=N(r,\psi)+S(r,f)\notag\\
&\leq\sum_{m=2}N_{(m,1)}(r,F_{1})-\sum_{m=2}\overline{N}_{(m,1)}(r,F_{1})+S(r,f),\label{eq14.13}
\end{align}
as $r\not\in E_{1}\cup E_{2}\cup E_{3}\cup E_{4}$ and $r\to\infty$. Then \eqref{eq14.12} and \eqref{eq14.13} deduce
\begin{align}
\overline{N}(r,F_{1})\leq N_{(1,1)}(r,F_{1})+S(r,f),\label{eq14.14}
\end{align}
as $r\not\in E_{1}\cup E_{2}\cup E_{3}\cup E_{4}$ and $r\to\infty$. Likely, we also have
\begin{align}
\overline{N}\left(r,\frac{1}{F_{1}-c}\right)\leq N_{(1,1)}\left(r,\frac{1}{F_{1}-c}\right)+S(r,f)\label{eq14.15}
\end{align}
and
\begin{align}
\overline{N}\left(r,\frac{1}{F_{1}-1}\right)&\leq N_{(1,1)}\left(r,\frac{1}{F_{1}-1}\right)+S(r,f),\label{eq14.16}
\end{align}
as $r\not\in E_{1}\cup E_{2}\cup E_{3}\cup E_{4}$ and $r\to\infty$. If $\varphi\equiv\psi$, then $F_{1}$ and $G_{1}$ share four values CM, and hence $f^{(k)}$ and $g^{(k)}$ share four values CM. Hence, we assume that $\varphi\not\equiv\psi$. So according to a simple computation, we can know by \eqref{eq14.6} and \eqref{eq14.7} that
\begin{align}
& N_{(1,1)}\left(r,\frac{1}{F_{1}}\right)+N_{(1,1)}\left(r,\frac{1}{F_{1}-c}\right)+N_{(1,1)}\left(r,\frac{1}{F_{1}-1}\right)\leq N\left(r,\frac{1}{\varphi-\psi}\right)\notag\\
&\leq T(r,\varphi-\psi)= m(r,\varphi-\psi)+N(r,\varphi-\psi)\leq N(r,\varphi-\psi)+S(r,f)\notag\\
&\leq \sum_{n=2}N_{(1,n)}(r,G_{1})-\sum_{n=2}\overline{N}_{(1,n)}(r,G_{1})+\sum_{m=2}N_{(m,1)}(r,F_{1})-\sum_{m=2}\overline{N}_{(m,1)}(r,F_{1})\notag\\
&+S(r,f)=N(r,F_{1})-\overline{N}(r,F_{1})+S(r,f),\label{eq14.17}
\end{align}
as $r\not\in E_{1}\cup E_{2}\cup E_{3}\cup E_{4}$ and $r\to\infty$. Thus,  we apply (ii) of Lemma 5 and \eqref{eq14.14}-\eqref{eq14.17} that
\begin{align}
&2T(r,F_{1})=\overline{N}(r,F_{1})+\overline{N}\left(r,\frac{1}{F_{1}}\right)+\overline{N}\left(r,\frac{1}{F_{1}-c}\right)+\overline{N}\left(r,\frac{1}{F_{1}-1}\right)+S(r,F)\notag\\
&\leq N_{(1,1)}\left(r,\frac{1}{F_{1}-1}\right)+N_{(1,1)}\left(r,\frac{1}{F_{1}-c}\right)+N_{(1,1)}\left(r,\frac{1}{F_{1}-d}\right)+\overline{N}(r,F_{1})\notag\\
&+S(r,f)\leq N(r,F_{1})+S(r,f)\leq T(r,F_{1})+S(r,f),\label{eq14.18}
\end{align}
as $r\not\in E_{1}\cup E_{2}\cup E_{3}\cup E_{4}$ and $r\to\infty$. Thus, we obtain from \eqref{eq14.18} that 
\begin{align}
T(r,F_{1})=S(r,f),\label{eq14.19}
\end{align}
as $r\not\in E_{1}\cup E_{2}\cup E_{3}\cup E_{4}$ and $r\to\infty$. This is a contradiction. Therefore,  we obtain that $F_{1}$ and $G_{1}$ share four values CM, and hence $f^{(k)}$ and $g^{(k)}$ share four values CM.

\par
\vskip 2mm From above discussions, we conclude that $f^{(k)}\equiv g^{(k)}$, or $f$ and $g$ are functions with normal growth, have the same order of growth, a positive integer or infinity.
\par
\vskip 2mm {\bf Case 2}\quad  $a_{1}, a_{2}$ are two distinct finite values, and $a_{3}$ is a non-constant small function.  Let 
$$F=f^{(k)}-a_{3},\quad G=g^{(k)}-a_{3},$$
and then it follows from above that $F$ and $G$ share $0,\infty, b=a_{1}-a_{3}, c=a_{2}-a_{3}$ IM*. 
By (2), (9), (i)-(ii) of Lemma 5 and Lemma 3, we have
\begin{eqnarray*}
\begin{aligned}
&2T(r,F)=\bar{m}_{0,q}(r,F)+\sum_{d\in\widehat{\mathbb{C}}}\widehat{N}(r,\frac{1}{F-d_{i}})+S(r,F)=2T(r,F)+\sum_{i=1}^{q-2}m(r,\frac{1}{F-d_{i}})\\
&-\overline{N}(r,F)-\overline{N}(r,\frac{1}{F})+\sum_{i=1}^{p-2}N(r,\frac{1}{F-d_{i}})-\sum_{i=1}^{p-2}\overline{N}(r,\frac{1}{F-d_{i}})+S(r,F)\\
&=2T(r,F)-2T(r,F)+\overline{N}(r,\frac{1}{F-b})+\sum_{i=1}^{q-2}m(r,\frac{1}{F-d_{i}})+\sum_{i=1}^{p-2}N(r,\frac{1}{F-d_{i}})\\
&+\overline{N}(r,\frac{1}{F-c})-(p-2-\varepsilon_{3})T(r,F)+\varepsilon_{1}T(r,F)\leq \overline{N}(r,\frac{1}{F-b})+\overline{N}(r,\frac{1}{F-c})\\
&+\varepsilon_{4} T(r,F)\leq \overline{N}(r,\frac{1}{F-b})+\overline{N}(r,\frac{1}{F-c})+\varepsilon_{4} T(r,F),
\end{aligned}
\end{eqnarray*}
which yields
\begin{align}
2T(r,F)\leq \overline{N}(r,\frac{1}{F-b})+\overline{N}(r,\frac{1}{F-c})+\varepsilon_{4} T(r,F)
\end{align}
for all $r>e$ outside a set $E_{2}\subset(e,\infty)$ of logarithmic density $0$. Similarly, we also have
\begin{align}
2T(r,G)\leq \overline{N}(r,\frac{1}{G-b})+\overline{N}(r,\frac{1}{G-c})+\varepsilon_{4} T(r,F)
\end{align}
for all $r>e$ outside a set $E_{3}\subset(e,\infty)$ of logarithmic density $0$. And hence
\begin{align}
\overline{N}(r,F)+\overline{N}(r,\frac{1}{F})+\overline{N}(r,\frac{1}{G})+\overline{N}(r,F)=\varepsilon_{4} T(r,F).
\end{align}
We suppose 
$$F_{1}=\frac{F-b}{c-b}, \quad G_{1}=\frac{G-b}{c-b},$$
and so the fact that $F$ and $G$ share $0,\infty, b=a_{1}-a_{3}, c=a_{2}-a_{3}$ IM* can imply $F_{1}$ and $G_{1}$ share $0,1,\infty$ and $D=\frac{-b}{c-b}\not\equiv 0,1,\infty$ IM*. And $D$ is a non-constant small function. Thus, we can easy to get from (44) that
\begin{align}
\overline{N}(r,F_{1})+\overline{N}(r,\frac{1}{F_{1}-D})+\overline{N}(r,\frac{1}{G_{1}-D})+\overline{N}(r,G_{1})=\varepsilon_{4} T(r,F).
\end{align}
By Lemma 7, we have that $F_{1}$ is a quasi-M$\ddot{o}$bius transformation of $G_{1}$. Then Lemma 8 can yield $F_{1}$ and $G_{1}$ share $0,1,\infty$ CM* and $D$ IM*, and hence by Lemma D1, we obtain $F_{1}$ and $G_{1}$ share four values CM. That is to say,  either $f^{(k)}\equiv g^{(k)}$, or $f$ and $g$ are functions with normal growth, have the same order of growth, a positive integer or infinity.
\par
\vskip 2mm {\bf Case 3} $a_{1}$ is a finite values, and $ a_{2},a_{3}$ are two distinct non-constant small functions. Then combining (2), (9), (i)-(ii) of Lemma 5 and Lemma 3, we have
\begin{eqnarray*}
\begin{aligned}
&2T(r,f^{(k)})=\bar{m}_{0,q}(r,f^{(k)})+\sum_{d\in\widehat{\mathbb{C}}}\widehat{N}(r,\frac{1}{f^{(k)}-d_{i}})+S(r,f^{(k)})=2T(r,f^{(k)})\\
&+\sum_{i=1}^{q-2}m(r,\frac{1}{f^{(k)}-d_{i}})+\sum_{i=1}^{p-2}N(r,\frac{1}{f^{(k)}-d_{i}})-\sum_{i=1}^{p-2}\overline{N}(r,\frac{1}{f^{(k)}-d_{i}})-\overline{N}(r,f^{(k)})\\
&-\overline{N}(r,\frac{1}{f^{(k)}-a_{1}})+S(r,f^{(k)})=2T(r,f^{(k)})+\sum_{i=1}^{q-2}m(r,\frac{1}{f^{(k)}-d_{i}})-\overline{N}(r,f^{(k)})\notag\\
&-\overline{N}(r,\frac{1}{f^{(k)}-a_{1}})+\sum_{i=1}^{p-2}N(r,\frac{1}{f^{(k)}-d_{i}})-(p-2-\varepsilon_{3})T(r,f^{(k)})+S(r,f^{(k)})\\
&\leq (2+\varepsilon_{3}) T(r,f^{(k)})-\overline{N}(r,f^{(k)})-\overline{N}(r,\frac{1}{f^{(k)}-a_{1}}),
\end{aligned}
\end{eqnarray*}
which yields
\begin{align}
\overline{N}(r,f^{(k)})+\overline{N}(r,\frac{1}{f^{(k)}-a_{1}})= \varepsilon_{3} T(r,f^{(k)}).
\end{align}

for all $r>e$ outside a set $E_{2}\subset(e,\infty)$ of logarithmic density $0$. Similar to prove as in {\bf Case 2}, we can get that either $f^{(k)}\equiv g^{(k)}$, or $f$ and $g$ are functions with normal growth, have the same order of growth, a positive integer or infinity.
\par
\vskip 2mm {\bf Case 4} $a_{1}, a_{2}, a_{3}$ are three distinct non-constant small functions. By (2), (9), (i)-(ii) of Lemma 5 and Lemma 3, we have
\begin{eqnarray*}
\begin{aligned}
2T(r,f^{(k)})&=\bar{m}_{0,q}(r,f^{(k)})+\sum_{d\in\widehat{\mathbb{C}}}\widehat{N}(r,\frac{1}{f^{(k)}-d_{i}})+S(r,f^{(k)})\\
&=T(r,f)+\sum_{i=1}^{q-1}m(r,\frac{1}{f^{(k)}-d_{i}})+\sum_{i=1}^{p-1}N(r,\frac{1}{f^{(k)}-d_{i}})\\
&-\overline{N}(r,f^{(k)})-\sum_{i=1}^{p-1}\overline{N}(r,\frac{1}{f^{(k)}-d_{i}})+S(r,f^{(k)})\\
&=T(r,f^{(k)})+\sum_{i=1}^{q-1}m(r,\frac{1}{f^{(k)}-d_{i}})-\overline{N}(r,f^{(k)})\notag\\
&+\sum_{i=1}^{p-1}N(r,\frac{1}{f^{(k)}-d_{i}})-(p-3-\varepsilon_{3})T(r,f^{(k)})+S(r,f^{(k)})\\
&\leq T(r,f^{(k)})+\sum_{i=1}^{q-3}m(r,\frac{1}{f^{(k)}-d_{i}})\\
&-\sum_{i=1}^{p-3}m(r,\frac{1}{f^{(k)}-d_{i}})-\overline{N}(r,f^{(k)})+\varepsilon_{3} T(r,f^{(k)})\\
&\leq T(r,f^{(k)})-\overline{N}(r,f^{(k)})+\varepsilon_{3} T(r,f^{(k)}),
\end{aligned}
\end{eqnarray*}
which yields
\begin{align}
T(r,f^{(k)})\leq \varepsilon_{3} T(r,f^{(k)})
\end{align}
for all $r>e$ outside a set $E_{2}\subset(e,\infty)$ of logarithmic density $0$. Take $\varepsilon=\frac{1}{12}$, we can get from (47) that $T(r,f^{(k)})=S(r,f^{(k)})$, but it is impossible.
\par
\vskip 2mm From above four cases discussions, we conclude that $f^{(k)}\equiv g^{(k)}$, or $f$ and $g$ are functions with normal growth, have the same order of growth, a positive integer or infinity. 
\par
\vskip 2mm This completes Theorem 8.

\

{\bf Acknowledgments} The author would like to thank anonymous referees for their helpful comments.

\
{\bf Data Availability Statement} Data sharing not applicable to this paper as no data sets were generated or analyzed during the current study.

\
{\bf Conflict of interest} The authors declare that they have no competing interests.


\end{document}